\documentclass{elsarticle}
\usepackage[utf8]{inputenc}
\usepackage[T1]{fontenc}
\usepackage[english]{babel}
\usepackage[style=english]{csquotes}
\usepackage{amsmath}
\usepackage{amsthm}
\usepackage{amsfonts,amssymb,amsopn,amstext,amscd,epsf,mathrsfs}
\usepackage{latexsym}
\usepackage{microtype}
\usepackage{enumitem}
\usepackage{tikz-cd}
\usepackage{tikz}
\usetikzlibrary{arrows,calc,positioning}
\usepackage{array}
\usepackage{lipsum}

\usepackage{natbib}

\numberwithin{equation}{section}
\theoremstyle{definition}
\newtheorem{defn}[equation]{Definition}{}
{}
\newtheorem{esempi*}[]{Example}{}
{}
\newtheorem{rmk}[equation]{Remark}{}

\newtheorem{thm}[equation]{Theorem}{}
\newtheorem*{thm*}{Theorem} 
\newtheorem*{lem*}{Lemma} 
\newtheorem{lem}[equation]{Lemma}
\newtheorem{prop}[equation]{Proposition}{}
\newtheorem{cor}[equation]{Corollary}
{}

\makeatletter
\newcommand*\rel@kern[1]{\kern#1\dimexpr\macc@kerna}
\newcommand*\widebar[1]{%
  \begingroup
  \def\mathaccent##1##2{%
    \rel@kern{0.8}%
    \overline{\rel@kern{-0.8}\macc@nucleus\rel@kern{0.2}}%
    \rel@kern{-0.2}%
  }%
  \macc@depth\@ne
  \let\math@bgroup\@empty \let\math@egroup\macc@set@skewchar
  \mathsurround\z@ \frozen@everymath{\mathgroup\macc@group\relax}%
  \macc@set@skewchar\relax
  \let\mathaccentV\macc@nested@a
  \macc@nested@a\relax111{#1}%
  \endgroup
}

\newcommand{\PP}{\mathbb{P}^1_k}
\newcommand{\Sp}{\textnormal{Spec}} 
\newcommand{\Rp}{\textnormal{Rep}} 
\newcommand{\St}{\textnormal{Strat}} 
\newcommand{\Sets}{\textnormal{Sets}} 
 
\newcommand{\Vc}{\textnormal{Vec}} 

\newcommand{\CC}{\mathcal{C}} 
 
\newcommand{\OO}{\mathcal{O}}
\newcommand{\FF}{\mathcal{F}}
\newcommand{\GG}{\mathcal{G}}
\newcommand{\EE}{\mathcal{E}}

\newcommand{\DD}{\mathcal{D}}
\newcommand{\hm}{\text{H}}

\newcommand{\ZZ}{\mathbb{Z}}
\newcommand{\Aut}{\textnormal{Aut}}
\usepackage[hidelinks]{hyperref}
\hypersetup{
  colorlinks   = true, 
  urlcolor     = black, 
  linkcolor    = black,
  citecolor   = black, 
}

\title{Specialization map between stratified bundles and pro-\'etale fundamental group}
\author{Elena Lavanda}

\begin{document}

\begin{abstract}
    Given a projective family of semi-stable curves over a complete discrete valuation ring of characteristic $p>0$  with algebraically closed residue field, we construct a specialization functor between the category of continuous representations of the pro-\'etale fundamental group of  the closed fibre and the category  of stratified bundles on the geometric generic fibre. 
By Tannakian duality, this functor induces a morphism between the corresponding affine group schemes. We show that this morphism is a lifting of the specialization map, constructed by Grothendieck, between the \'etale fundamental groups.

\


\end{abstract}

\maketitle

\section*{Introduction}

In \cite{Mumf}, given a complete discrete valuation ring $A$  of characteristic $p>0$  with fraction field $K$ and residue field $k$, Mumford associated with a flat Schottky group $G\subset \text{PGL}_2(K)$  a stable curve $X$ over $A$ with $k$-split degenerate closed fibre $X_0$ and non-singular generic fibre $X_K$,  such that $G$ is the group of covering transformations of the universal cover $Y_0$ of $X_0$. Moreover, he proved that every such curve $X$ can be constructed in this way for a unique flat Schottky group $G$ and that, if $X$ has arithmetic genus $g$,  $G$ is a free group with $g$ generators.

This setting was later used by Gieseker in \cite{Gies1} to prove that, for any prime $p>0$ and every integer $g>1$, there exists a stable curve of arithmetic genus $g$ in characteristic $p$ that admits a semi-stable bundle of rank two whose Frobenius pull-back is not semi-stable. Given a stable curve $X$ over $A$ with $k$-split degenerate closed fibre  and non-singular generic fibre, he introduced the notion of coherent sheaves with meromorphic descent data on the universal cover of the completion $\widehat{X}$ of $X$ along its closed fibre and he proved that the category they form is equivalent to the category of coherent sheaves on the generic fibre $X_K$. Then he associated with each $K$-linear representation of the group $G$, constructed by Mumford, a sheaf with meromorphic descent data, and hence, via the equivalence of categories, a bundle on the generic fibre. Furthermore, repeating the argument for all Frobenius twists of $X$, he associated with each representation of $G$ a stratified bundle on the geometric generic fibre. Finally, by cleverly choosing the representation, he  was able to construct a semi-stable bundle with the required properties. 

In this article we generalize Gieseker's  construction of stratified bundles from representations by removing the assumption on the degeneracy of the closed fibre.

 In the first section, we present an explicit computation of the pro-\'etale fundamental group of a projective normal crossing curve defined over an algebraically closed field.

\begin{thm*} \textnormal{(See \hyperref[genlemmapro]{Theorem. \ref*{genlemmapro}})} Given $X$ a connected projective normal crossing curve defined over an algebraically closed field and  $\xi$  a geometric point of $X$, let, for $j=1,\dots,N$,  $C_j$ be the irreducible components of $X$, $\widebar{C_j}$  their normalization and $\xi_j$ a fixed   geometric point  for every $\widebar{C_j}$, then
\[\pi^{\text{pro\'et}}_1(X,\xi)\simeq\mathbb{Z}^{\star |I|-N+1}\star_N \pi^{\text{\'et}}_1(\widebar{C_1},\xi_1)\star_N\cdots\star_N\pi^{\text{\'et}}_1(\widebar{C_N},\xi_N),\]
where $I$ is the set of singular points of $X$, $\mathbb{Z}^{\star |I|-N+1}$ is the free product of $|I|-N+1$ copies of $\mathbb{Z}$ and $\star_N$ is the co-product of Noohi groups.
\end{thm*}

In particular, if $X_0$ is a degenerate stable curve over an algebraically closed field $k$ of characteristic $p>0$, we see that its pro-\'etale fundamental group is isomorphic to the Schottky group  defined by Mumford and hence we have a more geometrical interpretation of the latter.

In the second section, given a topological group $G$, we illustrate the properties of its algebraic hull $G^{\text{cts}}$, which is defined as the affine group scheme associated with the Tannakian category of continuous representations of $G$.  In particular, we give an explicit description of the algebraic hull of a pro-finite group.

In the third section, given a complete DVR $A$ of equicharacteristic $p$ with algebraically closed residue field, we set $X$ to be a projective semi-stable curve over $\Sp(A)$ with connected closed fibre $X_0$ and smooth generic fibre $X_K$. We associate with a $K$-linear continuous representation $\rho$ of the pro-\'etale fundamental group of the closed fibre a geometric covering  $Y_{\rho}$ of the completion $\widehat{X}$ of $X$ along its closed fibre. Moreover, we show that meromorphic descent data on coherent sheaves over $\mathcal{Y}_{\rho}$ descend to coherent sheaves on $X_K$.

In  the fourth section, we extend this result to stratified bundles and this leads us to the definition of a specialization functor.

\begin{thm*} \textnormal{(See \hyperref[groupmor]{Theorem \ref*{groupmor}})} Let $X$ be a projective semi-stable curve over $\Sp(A)$ with connected closed fibre and smooth generic fibre, then the descent of stratified bundles with meromorphic descent data induces a tensor functor
\[\text{sp}_{\widebar K}\colon\Rp^{\textnormal{cts}}_{\widebar K}(\pi_1^{\text{pro\'et}}(X_0,\xi))\to\St(X_{\widebar K}),\]
which, by Tannakian duality, corresponds to a morphism of group schemes over $\widebar K$
\[\text{sp}\colon\pi^{\text{strat}}(X_{\widebar K})\to (\pi_1^{\text{pro\'et}}(X_0,\xi))^{\textnormal{cts}}.\]
\end{thm*}

We conclude by showing that this morphism of group schemes  is a lifting of the specialization map between the \'etale fundamental groups of $X_{\widebar K}$ and $X_0$ constructed by Grothendieck in \cite{SGA1}.

\

\noindent\textit{Acknowledgments.} The results contained in this article are part of my PhD thesis, written under the supervision of Hélène Esnault. I would like to thank her for guiding me through my graduate studies, for her support and for the many enlightening discussions we had. 

\section{Pro-\'etale fundamental group of semi-stable curves}\label{c1}

In \cite{BS} the authors introduced the notion of tame infinite Galois categories and proved that every such category is equivalent to the category of sets with a continuous action of a Noohi group. This notion generalizes the concept of Galois categories, introduced by Grothendieck in \cite{SGA1} and it is used to construct the pro-\'etale fundamental group of a scheme. Before defining this group, we recall the definition and basic properties of Noohi groups.

\begin{defn}Let $G$ be a topological group and $F_G\colon G\text{-}\Sets\to\Sets$ be the forgetful functor, we say that $G$ is a \textit{Noohi group} if the  natural map $G\to \Aut(F_G)$ is an isomorphism of topological groups, where $\Aut(F_G)$ is topologized  by the compact-open topology on $\Aut(S)$ for all $S\in\Sets$.  
\end{defn}

\begin{defn} Given $G$ a topological group, we define the \textit{Ra\u{\i}kov completion of $G$}, which is denoted by $\widehat{G}^R$, as its completion with respect to its two-sided uniformity  (see \cite{AT}). 
We say that a  topological group $G$ is \textit{Ra\u{\i}kov complete} if the natural morphism $\sigma\colon G\to\widehat{G}^R$, constructed in \cite[Thm. 3.6.10]{AT}, is an isomorphism.
\end{defn}

\begin{prop}[\cite{BS}, Prop. 7.1.5]\label{Noocomplete} Let $G$ be a topological group with a basis of open neighborhoods of $1\in G$ given by open subgroups and $F_G\colon G\text{-}\Sets\to\Sets$  the forgetful functor, then  $\Aut(F_G)$ is naturally isomorphic to $\widehat{G}^R$. Hence, $G$ is a Noohi group if and only if it is Ra\u{\i}kov complete.
\end{prop}

The pro-\'etale fundamental group of a scheme $X$, in analogy with the \'etale fundamental group, is defined as the Noohi group associated with the category of geometric coverings of $X$. 

\begin{defn} \label{geomcovdef} Given $X$ a  locally topologically  Noetherian connected scheme, we call \textit{geometric covering of X} any \'etale $X$-scheme $Y$ such the structure map $Y\to X$ satisfies the valuative criterion of properness. 
We denote by $\text{Cov}_X$ the category of geometric coverings, where the maps are given by $X$-morphisms.
\end{defn}

\begin{thm}[\cite{BS}, Lemma 7.4.1]\label{proet} Let $X$ be a locally topologically Noetherian connected scheme, $\xi$ a geometric point of $X$ and set $\textnormal{ev}_{\xi}$ to be the following functor
\[\textnormal{ev}_{\xi}\colon \textnormal{Cov}_X\to\Sets\text{, }\textnormal{ev}_{\xi}(\pi\colon Y\to X)=\pi^{-1}(\xi),\]
then  the group $\textnormal{Aut}(\textnormal{ev}_{\xi})$, endowed with the compact-open topology, is a Noohi group.
Moreover, the functor $\textnormal{ev}_{\xi}$ induces an equivalence of categories 
\[\textnormal{ev}_{\xi}\colon\textnormal{Cov}_X\simeq\textnormal{Aut}(\textnormal{ev}_{\xi})\textnormal{-}\Sets.\]
\end{thm}

\begin{defn} Given $X$ a locally topologically Noetherian connected scheme and  $\xi$  a geometric point of $X$, we define the \textit{pro-\'etale fundamental group} of $X$, as in \cite[Def. 7.4.2]{BS}, to be the group
\[\pi_1^{\text{pro\'et}}(X,\xi):=\text{Aut}(\text{ev}_{\xi}).\] 
\end{defn}

From the pro-\'etale fundamental group, we can retrieve both the enlarged fundamental defined in \cite{SGA3} and the \'etale fundamental group.

\begin{prop}[\cite{BS}, Lemma 7.4.3 and Lemma 7.4.6]\label{complet} Let $X$ be a locally topologically Noetherian connected scheme and  $\xi$  a geometric point of $X$, then
\begin{itemize}
\item the pro-discrete completion of $\pi_1^{\text{pro\'et}}(X,\xi)$ is isomorphic to the enlarged fundamental group $\pi_1^{\text{SGA3}}(X,\xi)$,
\item the pro-finite completion of $\pi_1^{\text{pro\'et}}(X,\xi)$ is isomorphic to the \'etale fundamental group $\pi_1^{\text{\'et}}(X,\xi)$.
\end{itemize}
\end{prop}

\begin{prop}[\cite{BS}, Lemma 7.4.10]\label{geomuni} If $X$ is geometrically unibranch, then 
\[\pi_1^{\text{pro\'et}}(X,\xi)\simeq\pi_1^{\text{\'et}}(X,\xi).\]
\end{prop}

Before computing the pro-\'etale fundamental group of normal crossing curve, we state some basic definitions. 

\begin{defn} Let $C$ be a scheme of dimension $1$ of finite type over an algebraically closed field $\widebar F$, then $C$ is a \textit{semi-stable curve} if it is reduced  and its singular points are ordinary double points. If $F$ is any field and  $\widebar F$ is a fixed algebraic closure of $F$, then a curve $C$ over $F$ is called \textit{semi-stable}  if $C_{\widebar F}=C\times_F \Sp(\widebar F)$ is a semi-stable  curve over $\widebar F$.
\end{defn}

\begin{defn}\label{normcrossalgcl} Let $C$ be a scheme of dimension $1$ of finite type over an algebraically closed field $\widebar F$, then we say that $C$ is  a \textit{normal crossing curve} if its associated reduced scheme $C_{\text{red}}$ is a semi-stable curve. If $F$ is any field and  $\widebar F$ is a fixed algebraic closure of $F$, then a curve $C$ over $F$ is called \textit{normal crossing}  if its base change $C_{\widebar F}$ is a normal crossing curve over $\widebar F$.
\end{defn}

\begin{defn} Given a scheme $S$,  a \textit{semi-stable curve over $S$} is a flat scheme $X$  over $S$, whose fibres are semi-stable curves. 
\end{defn}

The main idea behind the  computation  of the pro-\'etale fundamental group of normal crossing curves is to generalize  \cite[Exp. IX Cor. 5.4]{SGA1} in terms of the pro-\'etale fundamental group. Hence, we need an explicit construction of the co-product of Noohi groups. 

\begin{rmk}[\cite{BS}, Example 7.2.6]\label{deffreeprodNoohi} Given two Noohi groups $G$ and $H$, we  set $\CC_{G,H}$ to be the category of triples $(S, \rho_G, \rho_H)$ where $S$ is a set and $\rho_G$, $\rho_H$ are continuous actions on  $S$ of $G$ and $H$ respectively. Let $\text{forg}\colon\CC_{G,H}\to\Sets$ be the forgetful functor, then the group $\text{Aut}(\text{forg})$ is Noohi and it is, in fact, the co-product of $G$ and $H$ in the category of Noohi groups. We will denote the co-product of Noohi groups $G$ and $H$ by $G\star_N H:=\text{Aut}(\text{forg})$.
\end{rmk}

We  give now an alternative description of the co-product in the category of Noohi groups. In what follows, given two topological groups $G$ and $H$, we denote by  $G\star H$ the co-product in the category of topological groups,  constructed in \cite{graev}.

\begin{lem}\label{basisneigh} For two Noohi groups $G$ and $H$ with a basis of open neighborhoods of $1$ given by open subgroups, we set $\mathcal{B}$ to be the collection of open subsets of   $G\star H$ of the form
\[x_1\Gamma_1 y_1\cap\cdots\cap x_n\Gamma_n y_n,\]
with $n\in\mathbb{N}$, $x_i,y_i\in G\star H$ and $\Gamma_i\subseteq G\star H$ open subgroups of $G\star H$. If we restrict the topology on $G\star H$ to the topology induced by $\mathcal{B}$, we obtain a topological group $G\star_{\mathcal{B}} H$ with a basis of open neighborhoods of $1\in G\star H$ given by open subgroups. 
\end{lem}
\proof  Given $x,y\in G\star H$ and $\Gamma\subset G\star H$ an open subgroup,let $m$ be the group operation, then $(z_1,z_2)\in m^{-1}(x\Gamma y)$ implies that
\[yz_2^{-1}z_1^{-1}x=(x^{-1}z_1z_2y^{-1})^{-1}\in \Gamma.\]
Hence, the multiplication is continuous because we have, for every $x,y$ and $\Gamma$,
\[(z_1,z_2)\in x\Gamma y z_2^{-1}\times z_1^{-1}x\Gamma y\subset m^{-1}(x\Gamma y).\]
Let $i$ be the inverse morphism, then $ y^{-1}\Gamma x^{-1}\subset i^{-1}(x\Gamma y)$, for every $x,y$ and every $\Gamma$, thus $G\star_{\mathcal{B}} H$ is a topological group.

To conclude, it suffices to show that every set $x\Gamma y\in \mathcal{B}$ such that $1\in x\Gamma y$ contains an open subgroup of $G\star_{\mathcal{B}} H$. The condition $1\in x\Gamma y$ implies that $x^{-1}y^{-1}\in\Gamma$. 
The set $y^{-1}\Gamma y$ is, by definition, an open subgroup of  $G\star_{\mathcal{B}} H$. Moreover, we see that  $y^{-1}\Gamma y\subset x\Gamma y$  because, given $\delta\in y^{-1}\Gamma y$, we have,  for some $\gamma\in \Gamma$,
\[\delta=y^{-1}\gamma y=x(x^{-1}y^{-1})\gamma y\in  x\Gamma y.\] 
\endproof

\begin{cor}\label{noohicoproductcompl} Let $G$ and $H$ be two Noohi groups with a basis of open neighborhoods of $1$ given by open subgroups, then the co-product in the category of Noohi groups $G\star_N H$ is isomorphic to the Ra\u{\i}kov completion of the topological group $G\star_{\mathcal{B}} H$, defined above.
\end{cor}
\proof By \hyperref[basisneigh]{Lemma \ref*{basisneigh}},  ${G\star_{\mathcal{B}} H}$ has  a basis of open neighbourhoods of $1$ given by
open subgroups. Hence, by \hyperref[Noocomplete]{Proposition \ref*{Noocomplete}}, it suffices to prove that the categories $G\star_N H$-$\Sets$ and ${G\star_{\mathcal{B}} H}$-$\Sets$ are equivalent.

By the universal property of the co-product of topological groups, the categories $G\star_N H$-$\Sets$ and $G\star H$-$\Sets$ are equivalent. 
Furthermore, the identity induces a continuous morphism $G\star H\to G\star_{\mathcal{B}} H$, which corresponds to a fully faithful functor
$G\star_{\mathcal{B}} H\text{-}\Sets\to G\star H\text{-}\Sets$.
Let $\rho$ be a continuous action of $G\star H$ on a set $S$, then the map $\rho\colon G\star H\to \text{Aut}(S)$ is continuous with respect to the compact-open topology on $\text{Aut}(S)$. Since a basis of open neighborhoods of $1\in\text{Aut}(S)$ is given by stabilizers of finite subsets of $S$,  the inverse image via $\rho$ of any open neighborhood of $1\in\text{Aut}(S)$ contains an open subgroup of $G\star H$. By construction, this implies that the map $\rho$ is continuous also with respect to the topology of $G\star_{\mathcal{B}} H$, hence  the functor induced by the identity is an equivalence of category. 
\endproof

Note that, by the universal property, the co-product of two discrete groups in the category of topological groups is their abstract free product endowed with the discrete topology and it coincides with the co-product in the category of Noohi groups.

We proceed now with the computation of the pro-\'etale fundamental group of a normal crossing curve.

\begin{lem}\label{prored} Let $X$ be a locally Noetherian connected scheme, $X_{\text{red}}$ its associated reduced subscheme and  $\xi$  a geometric point of $X$, then
\[\pi_1^{\text{pro\'et}}(X_{\text{red}},\xi)\simeq\pi_1^{\text{pro\'et}}(X,\xi).\]
\end{lem}
\proof By \cite[Thm. 18.1.2]{EGA4} the category of schemes that are \'etale over $X$ is equivalent to the category of schemes that are \'etale over  $X_{\text{red}}$. Thus, it suffices to prove that an \'etale scheme $Y$ over $X$ satisfies the valuative criterion of properness if and only if  $Y\times_X X_{\text{red}}=Y_{\text{red}}$ does.

Let $R$ be any discrete valuation with fraction field $F$, then any morphism $\Sp(F)\to Y$ factors through $Y_{\text{red}}$ and similarly any morphism $\Sp(R)\to X$ factors through $X_{\text{red}}$. Hence, it is clear that, for any diagram of the form
\[
\begin{CD}
\Sp(F) @>>> Y\\
@VVV @VVV \\
\Sp(R) @>>> X \text{,}\\
\end{CD}
\]
there exists a unique map $\Sp(A)\to Y$ that makes the diagram commutative if and only if there exist a unique map $\Sp(A)\to Y_{\text{red}}$ that makes the diagram between the associated reduced schemes commutative. 
\endproof

\begin{prop} \label{gendes} Let $g\colon X'\to X$ be a proper surjective morphism of finite presentation, then $g$ is a morphism of effective descent for geometric coverings. 
\end{prop}
\proof By \cite[Thm. 5.19]{Rydh} and \cite[Thm. 5.4]{Rydh}, $g$ is a morphism of effective descent for \'etale separated schemes.  Since geometric coverings are \'etale and satisfy the valuative criterion of properness, they are, in particular, separated \'etale morphisms.  Let $Y'$ be a geometric covering  of $X'$ with descent data relative to $g$, then $Y'$ descends to a separated \'etale  $X$-scheme $Y$. Moreover, since $g$ is proper, $Y'$ satisfies the valuative criterion of properness if and only if $Y$ does. Hence, $g$ is a morphism of effective descent for geometric coverings.
\endproof

\begin{prop} \label{genlemmapro} Let $X$ be a  projective connected normal crossing curve over an algebraically closed field $F$ and $\xi$   a geometric point of $X$. For $j=1,\dots N$, let  $C_j$ be the irreducible components of $X$, $\widebar{C_j}$  their normalizations and   $\xi_j$ a  fixed geometric point for every $\widebar{C_j}$, then
\[\pi^{\text{pro\'et}}_1(X,\xi)\simeq\mathbb{Z}^{\star |I|-N+1}\star_N \pi^{\text{\'et}}_1(\widebar{C_1},\xi_1)\star_N\cdots\star_N\pi^{\text{\'et}}_1(\widebar{C_N},\xi_N),\]
where $I$ is the set of singular points of $X$  and $\mathbb{Z}^{\star |I|-N+1}$ is the free product of $|I|-N+1$ copies of $\mathbb{Z}$.
\end{prop}
\proof By \hyperref[prored]{Lemma \ref*{prored}}, we can assume  that $X$ is a projective connected  semi-stable curve. Hence, by \hyperref[gendes]{Proposition \ref*{gendes}}, the normalization is a morphism of effective descent for geometric coverings. 
We prove the statement  by induction on $N$, the number of irreducible components of $X$. 

If $X$ is irreducible, the normalization $\widebar{X}$ is connected. In this simple setting the descent data of geometric coverings of $\widebar{X}$ with respect to the normalization can be described explicitly.
We denote by $(a_i,b_i)$ the pair of points of $\widebar{X}$ that are identified to $x_i\in I$ in $X$ and we set $F_{a_i}$ and $F_{b_i}$ to be the functors associating to each geometric covering its fibers over $a_i$ and $b_i$ respectively.
Giving descent data for $Y$, a geometric covering of $\widebar{X}$, with respect to the normalization is equivalent to giving a collection of bijections $\{\alpha_i\colon F_{a_i}(Y)\to F_{b_i}(Y)\}_{x_i\in I}$.

Let $\mathcal{C}$ be the category whose objects are given by the datum $(Y,\alpha_1,\dots,\alpha_r)$ with $Y$ a geometric covering of $\widebar{X}$ and $\alpha_i\colon F_{a_i}(Y)\to F_{b_i}(Y)$ isomorphisms of sets, and whose morphisms from $(Y,\alpha_i)$ to $(Z,\beta_i)$ are given by $\widebar{X}$-scheme morphisms $\varphi\colon Y\to Z$ such that, for every $i\in I$, the following diagram commutes
\[
\begin{CD}
F_{a_i}(Y) @>\alpha_i>>F_{b_i}(Y)\\
@VF_{a_i}(\varphi) VV @VVF_{b_i}(\varphi)V \\
F_{a_i}(Z) @>\beta_i>>F_{b_i}(Z)\ \text{.}\\
\end{CD}
\]
By construction, the category $\mathcal{C}$ is equivalent to the category of  geometric coverings of $X$. 
We claim that there exists an equivalence between the category $\mathcal{C}$ and the category $\mathcal{C}_{\mathbb{Z}^{\star r},\pi^{\text{pro\'et}}_1(\widebar{X},\xi_1)}$, defined as in \hyperref[deffreeprodNoohi]{Remark \ref*{deffreeprodNoohi}}, which is compatible with their fiber functors.
If the claim is true, then it follows that 
\[\pi^{\text{pro\'et}}_1({X},\xi)\simeq\mathbb{Z}^{\star |I|}\star_N\pi^{\text{pro\'et}}_1(\widebar{X},\xi_1).\] 

By definition,  $\pi^{\text{pro\'et}}_1(\widebar{X},\xi_1)=\text{Aut}(F_{\xi_1})$ acts on ${F}_{\xi_1}(Y)$ for every $(Y,\alpha_i)\in\CC$. 
Since $\widebar{X}$ is connected, we can choose, for every $i$, a path $\tau_i$ from $a_i$ to $b_i$ and a path $\sigma_i$ from $\xi_1$ to $a_i$ and we notice that every $\alpha_i\in\text{Hom}(F_{a_i}(Y),F_{b_i}(Y))$ can be written as 
\[\alpha_i=\tau_i\circ g_i \text{ for some }g_i\in\text{Aut}(F_{a_i}(Y)).\]
Hence, we can define the action $\rho_i$ of  $i$-th copy of $\mathbb{Z}$   on $F_{\xi_1}(Y)$ as 
\[\rho_i(1)=\sigma_i \circ g_i \circ \sigma_i^{-1},\]
which induces a functor 
\[\widetilde{F}_{\xi_1}(Y)\colon \CC\to\mathcal{C}_{\mathbb{Z}^{\star r},\pi^{\text{pro\'et}}_1(\widebar{X},\xi_1)}.\]

Given an object $(S,\rho_{1,\dots,r},\rho_{\xi_1})\in\mathcal{C}_{\mathbb{Z}^{\star r},\pi^{\text{pro\'et}}_1(\widebar{X},\xi_1)}$, there exists a geometric covering $Y$ of $\widebar{X}$ such that $F_{\xi_1}(Y)\simeq (S, \rho_{\xi_1}).$ 
Thus, we can define the following functor: 
\[G_{\xi_1}(S,\rho_i, \rho_{\xi_1})=(Y,\tau_i\circ \sigma_i^{-1}\circ \rho_i(1)\circ\sigma_i),\]
which clearly is a quasi-inverse functor of  $\widetilde{F}_{\xi_1}$.

Since, by construction, $\text{forg}\circ\widetilde{F}_{\xi_1}(Y)=F_{\xi_1}(Y)$, we have proved the previous claim.

Let us prove now the inductive step.
We fix $C_1$ an irreducible component of $X$ such that the geometric point $\xi$ does not lie in $C_1$ and such that $X\setminus C_1$ is connected. We denote by $I_1$  the set of pairs $(a^1_i,b^1_i)$  of points of $\widebar{C_1}$ identified to a singular point $x_i^1$ of $C_1$, then by the base case we conclude that 
\[\pi^{\text{pro\'et}}_1(C_1,\xi_1)\simeq \mathbb{Z}^{\star |I_1|}\star_N\pi^{\text{pro\'et}}_1(\widebar{C_1},\xi_1).\] 

We denote by $\widebar{X}_{N-1}$ be the complement of $\widebar{C_1}$ in the normalization of $X$, denoted by $\widebar X$, by $I_{N-1}$ the set of pairs $(a^{N-1}_i,b^{N-1}_i)$  of points of $\widebar{X}_{N-1}$ identified to a singular point $x_i^{N-1}$ of $X$ and  we set $X_{N-1}$ to be the curve obtained from $\widebar{X}_{N-1}$ identifying the  pairs in $I_{N-1}$. By construction, $X_{N-1}$ is a  projective connected semi-stable curve with $N-1$ irreducible components and, by the inductive hypothesis,
\[\pi^{\text{pro\'et}}_1(X_{N-1},\xi)\simeq\mathbb{Z}^{\star |I_{N-1}|-N+2}\star_N\pi^{\text{pro\'et}}_1(\widebar{C_2},\xi_2)\star_N\cdots\star_N \pi^{\text{pro\'et}}_1(\widebar{C_N},\xi_N).\]

Finally, we denote by $I_{1,N-1}$ the set of pairs $(a_i^{1},b_i^{N-1})$, with $a_i^1$  a point of $\widebar{C_1}$ and $b_i^{N-1}$ a point of $\widebar{X_{N-1}}$, that are identified in the remaining singular points of $X$. 
We fix a pair $(a_0^1,b_0^{N-1})\in I_{1,N-1}$ and we set $X'$ to be the curve obtained from gluing $C_1$ and $X_{N-1}$ along the pair $(a_0^1,b_0^{N-1})\in I_{1,N-1}.$

We define  $\mathcal{C}_0$ to be the category whose objects are triples $(Y_1,Y_{N-1},\alpha_0)$ with
 $Y_1$ a finite \'etale cover of $C_1$,  $Y_{N-1}$ a finite \'etale cover of $X_{N-1}$,
and $\alpha_0$ an isomorphism of sets $F_{a_0^1}(Y_1)\to F_{b_0^{N-1}}(Y_{N-1})$,
and  whose morphisms from $(Y_1,Y_{N-1},\alpha_0)$ to $(Z_1,Z_{N-1},\beta_0)$
are given by  pairs $(\varphi_1,\varphi_{N-1})$ with
 $\varphi_1\colon Y_1\to Z_1$   a morphism of $C_1$-schemes and
$\varphi_{N-1}\colon Y_{N_1}\to Z_{N-1}$  a morphisms $X_{N-1}$-schemes
 such that the following diagram commutes
\[
\begin{CD}
F_{a_0^1}(Y_1) @>\alpha_0>>F_{b_0^{N-1}}(Y_{N-1})\\
@VF_{a_0}(\varphi_1) VV @VVF_{b_0}(\varphi_{N-1})V \\
F_{a_0^1}(Z_1) @>\beta_0>>F_{b_0^{N-1}}(Z_{N-1})\ \text{.}\\
\end{CD}
\]
Clearly $\mathcal{C}_0$ is equivalent to the category of geometric coverings of $X'$ and we claim that the  categories  $\mathcal{C}_0$  and $\pi^{\text{pro\'et}}_1(C_1,\xi_1)\star_N\pi^{\text{pro\'et}}_1(X_{N-1},\xi)\text{-}\Sets$ are equivalent.

The group $\pi^{\text{pro\'et}}_1(C_1,\xi_1)=\Aut(F_{\xi_1})$  acts naturally on $F_{\xi_1}(Y_1)$.
Furthermore, the schemes $C_1$ and $X_{N-1}$ are connected, so we can choose the paths 
\[\sigma_1\colon F_{a_0^1}\to F_{\xi_1}\text{ and }\sigma_{N-1}\colon F_{b_0^{N-1}}\to F_{\xi}.\]
We call $\rho$ the action of $\pi^{\text{pro\'et}}_1(X_{N-1},\xi)\simeq\Aut(F_{\xi})$ on $F_{\xi}(Y_{N-1})$ and we  define, for every $g\in\Aut(F_{\xi})$,
\[\tau(g)=(\sigma_{N-1}\circ\alpha_0\circ\sigma_1^{-1})^{-1}\circ\rho(g)\circ(\sigma_{N-1}\circ\alpha_0\circ\sigma_1^{-1}).\] 
Then $\tau$ is  an action of $\Aut(F_{\xi})$ on $F_{\xi_1}(Y_1)$ and it induces a functor 
\[\widetilde{F}_{\xi_1}\colon\CC_0\to\pi^{\text{pro\'et}}_1(C_1,\xi_1)\star_N\pi^{\text{pro\'et}}_1(X_{N-1},\xi)\text{-}\Sets.\]

Given $(S,\rho_1,\rho_{N-1})\in\pi^{\text{pro\'et}}_1(C_1,\xi_1)\star\pi^{\text{pro\'et}}_1(X_{N-1},\xi)\text{-}\Sets$,  there exists  a geometric covering $Y_1$ of $C_1$ such that $F_{\xi_1}(Y_1)\simeq(S,\rho_1)$ and a geometric covering $Y_{N-1}$ of $X_{N-1}$ such that 
$F_{\xi}(Y_{N-1})\simeq (S,\rho_{N-1}).$
Thus, we can define the functor 
\[G_{\xi_1}(S,\rho_1,\rho_{N-1})=(Y_1,Y_{N-1},\sigma^{-1}_{N-1}\circ\text{Id}_S\circ \sigma_{1}),\]
which is a quasi-inverse of $\widetilde{F}_{\xi_1}$.

Finally, we observe that a geometric covering of $X$ corresponds to the datum of a geometric covering $Y$ of $X'$ and the isomorphisms $\alpha_i\colon F_{a_i^1}(Y)\to F_{b_i^{N-1}}(Y)$ for every remaining pair of points $\{a_i^1,b_i^{N-1}\}\in I_{1,N-1}$.
By the same argument of the base step, 
\[\pi_1^{\text{pro\'et}}(X,\xi)\simeq\mathbb{Z}^{\star |I_{1,N-1}|-1}\star_N\pi^{\text{pro\'et}}_1(X',\xi).\]
Hence, we obtain that
\[\pi_1^{\text{pro\'et}}(X,\xi)\simeq\mathbb{Z}^{\star|I|-N+1}\star_N\pi^{\text{pro\'et}}_1(\widebar{C_1},\xi_1)\star_N \cdots\star_N \pi^{\text{pro\'et}}_1(\widebar{C_N},\xi_N).\]
The statement follows because, since $\widebar{C_j}$ are normal, by  \hyperref[geomuni]{Proposition \ref*{geomuni}}, 
\[\pi^{\text{pro\'et}}_1(\widebar{C_j},\xi_j)\simeq \pi^{\text{\'et}}_1(\widebar{C_j},\xi_j).\]
 
\endproof

For the following sections, we consider a fixed isomorphism between $\pi_1^{\text{pro\'et}}(X,\xi)$ and $\mathbb{Z}^{\star |I|-N+1}\star_N \pi^{\text{\'et}}_1(\widebar{C_1},\xi_1)\star_N\cdots\star_N\pi^{\text{\'et}}_1(\widebar{C_N},\xi_N)$, as constructed in the previous proposition.

\section{Algebraic hulls}

\begin{defn}\label{ctsrep} Given a field $F$ and  a topological group $G$, a \textit{continuous  $F$-linear representation of }$G$ is a pair  $(V,\rho)$ of  a finite dimensional $F$-vector space  $V$ and an $F$-linear action $\rho\colon G\times V\to V$ that is continuous with respect to the discrete topology on $V$. We denote by $\Rp^{\text{cts}}_F(G)$ the category of continuous $F$-linear representations of $G$.
\end{defn}

\begin{defn} \label{topalghull} Let $F$ be a field and $G$ a topological group, we define the \textit{algebraic hull of $G$ over $F$} to be the affine group scheme over $F$ associated, by Tannakian duality (\cite[Thm. 2.11]{Miln}), with the pair $(\Rp^{\text{cts}}_F(G),\omega_G)$, where $\omega_G(V,\rho)=V$ is the forgetful functor. We denote this group scheme by $G^{\text{cts}}$.
\end{defn}

We recall the following  elementary result in topology theory.

\begin{lem}\label{repdiscrete} Given $F$ a field, $G$ a topological group, $V$ a finite dimensional $F$-vector space and $\rho\colon G\times V\to V$ an $F$-linear $G$-action on $V$, the following are equivalent:
\begin{enumerate}
\item $\rho\colon G\times V\to V$ is continuous with respect to the discrete topology on $V$,
\item the group morphism ${\rho}\colon G\to \Aut(V)$ is  continuous  with respect to the compact-open topology on $\Aut(V)$.
Moreover, the compact-open topology on $\Aut(V)$ coincides with the discrete topology on $\Aut(V)$.
\end{enumerate}
\end{lem}

\begin{rmk} In particular, if $G$ is a given topological group, $\widehat{G}^D$ is its  pro-discrete completion and $F$ any field, then there exists an equivalence of categories 
\[\Rp^{\text{cts}}_F(\widehat{G}^D) \to\Rp^{\text{cts}}_F(G).\]

Let $X$ be a locally topologically Noetherian connected scheme and $\xi$ a geometric point of $X$, by \hyperref[complet]{Proposition \ref*{complet}}, the pro-discrete completion of  $\pi^{\text{pro\'et}}_1(X,\xi)$ is isomorphic to $\pi^{\text{SGA3}}_1(X,\xi)$, hence it follows that, for every field $F$,
\[\Rp^{\text{cts}}_F(\pi^{\text{pro\'et}}_1(X,\xi)) \simeq \Rp^{\text{cts}}_F(\pi^{\text{SGA3}}_1(X,\xi)).\]
Note that this equivalence of categories holds even in the cases, presented for example in  \cite[Example 7.4.9]{BS}, where $\pi^{\text{pro\'et}}_1(X,\xi)$ and $\pi^{\text{SGA3}}_1(X,\xi)$ are not isomorphic as topological groups.
\end{rmk}

In the next statements we will describe  the algebraic hulls of finite and pro-finite groups. 

\begin{lem}\label{finitealg} Let $G$ be a finite group and $G^{\textnormal{cts}}$ be its algebraic hull over a given field $F$, then $G^{\textnormal{cts}}$ is isomorphic to the constant group scheme over $F$ associated with $G$.
\end{lem}
\proof  Since $G$ is finite, the category $\Rp^{\text{cts}}_F(G)$  is equivalent to the category of finite dimensional $F$-linear representations $\Rp_F(G)$ and hence to the category of finitely generated $F[G]$-modules, where $F[G]$ is the $F$-Hopf algebra generated by the elements of $G$. Let $F^G$ be the dual $F$-Hopf algebra of $F[G]$, then $\Rp_F(G)$ is equivalent to the category of finitely generated $F^G$-comodules. This implies, by \cite[Ex. 2.15]{Miln}, that $G^{\text{cts}}=\Sp(F^G)$ and hence it is the constant group scheme associated with $G$.
\endproof

\begin{lem}\label{commlim} Let $F$ be a field and $\pi=\varprojlim_{i}\pi_i$ be a complete pro-finite  group with surjective transition maps, then $\pi^{\textnormal{cts}}$, the algebraic hull of $\pi$ over $F$, is isomorphic to $F$- group scheme
\[\pi_F:=\varprojlim_i (\pi_i)_F,\]
where $(\pi_i)_F$ are the constant group schemes  associated with the finite quotients $\pi_i$.
\end{lem}
\proof Since $\pi_i$ is finite, by \hyperref[finitealg]{Lemma \ref*{finitealg}}, $\pi_i^{\text{cts}}$ is the constant group scheme over $F$ associated with $\pi_i$, which we denote by $(\pi_i)_F$. 

The natural map $\text{pr}_i\colon\pi\to\pi_i$ induces a tensor functor between the categories of continuous representations
\[F_{\varphi_i}\colon\Rp_F^{\text{cts}}(\pi_i)\to\Rp_F^{\text{cts}}(\pi)\text{, }F_{\varphi_i}(V,\rho):=(V,\rho\circ \text{pr}_i).\]
By \cite[Cor. 2.9]{Miln}, this induces, for each $i$, a morphism of $F$-group schemes
\[\varphi_i\colon\pi^{\text{cts}}\to (\pi_i)_F.\]
Hence, there exists a natural morphism of $F$-group schemes
\[\varphi\colon\pi^{\text{cts}}\to\varprojlim_i (\pi_i)_F,\]
which corresponds to a functor
\[F_{\varphi}\colon\Rp_F(\varprojlim_i (\pi_i)_F)\to\Rp_F(\pi^{\text{cts}})\simeq\Rp_F^{\text{cts}}(\pi).\]

Since the maps $\text{pr}_i$ are surjective, the functor $F_{\varphi_i}$ satisfies the criterion of \cite[Prop. 2.21]{Miln}. This implies that   $\varphi_i$  is faithfully flat and, in particular, that it is surjective.
If $\pi^{\text{cts}}=\Sp(A)$ and $(\pi_i)_F=\Sp(B_i)$, then the affine morphism $\varphi_i$ corresponds to an injective morphism of $F$-Hopf algebras $\varphi_i\colon B_i\subset A$. Thus, the induced map $\varinjlim_i B_i\to A$, which corresponds to the morphism $\varphi$, is injective as well and, by \cite[VI,Thm 11.1]{milneAGS}, 
is faithfully flat. Then, by \cite[Prop. 2.21.(a)]{Miln}, $F_{\varphi}$ is fully faithful and it remains to show that it is  also essentially surjective.

By \hyperref[repdiscrete]{Lemma \ref*{repdiscrete}}, given an object $(V,\rho)\in\Rp_F^{\text{cts}}(\pi)$, the map ${\rho}\colon\pi\to\Aut(V)$ is continuous with respect to the discrete topology on $\Aut(V)$, which implies that ${\rho}$ factors through a finite quotient $\pi_i$ of $\pi$. Thus, there exists a $(\pi_i)_F$-action $\rho_i$ on $V$, such that $\rho_i\circ \varphi_i=\rho.$
Let $ p_i\colon \varprojlim_i\pi_i^{\text{cts}}\to\pi_i^{\text{cts}}$ be the natural morphism of $F$- group schemes, then
\[F_{\varphi}(V,\rho_i\circ p_i):=(V,\rho_i\circ p_i\circ \varphi)=(V,\rho_i\circ \varphi_i)=(V,\rho).\]
\endproof

\section{Descent of coherent sheaves with meromorphic data}\label{c3}

The following notation will be used throughout  these last three sections.
We fix $k$ an algebraically closed field of characteristic $p>0$, we set $A$ to be a complete discrete valuation ring of characteristic $p$ with residue field $k$, we denote by $K$ the fraction field of $A$ and we set $S=\Sp(A)$. Moreover, we fix $X\to S$ a projective semi-stable curve with connected closed fibre $X_0$ and smooth generic fibre $X_K$.

\

Under the assumption that the closed fibre $X_0$ is degenerate, that is that the normalizations of its irreducible components are isomorphic to $\PP$, in \cite{Gies1} Gieseker associated with a $K$-linear representation of the free group $\ZZ^{\star r}$, with $r=p_a(X_0)$ the arithmetic genus of $X_0$, a stratified bundle on $X_{\widebar K}$. We have proved in \hyperref[genlemmapro]{Theorem \ref*{genlemmapro}} that the group $\ZZ^{\star r}$, which had only a computational description in \cite{Mumf}, is, in fact, isomorphic to the pro-\'etale fundamental group of $X_0$.

The degeneracy assumption was essential for Mumford because it allowed him to construct a universal cover of $X_0$. We can reinterpret this phenomenon also in terms of the pro-\'etale fundamental group. Indeed, if $X_0$ is degenerate, then the left regular $\pi_1^{\text{pro\'et}}(X_0,\xi)$-action on the set $S=\pi_1^{\text{pro\'et}}(X_0,\xi)$ is continuous with respect to the discrete topology on $S$. Hence, it induces an object of the category $\pi_1^{\text{pro\'et}}(X_0,\xi)$-$\Sets$, which corresponds to the universal cover $Y_0$ of $X_0$.
On the other hand, if $X_0$ is not degenerate, the regular action on $S=\pi_1^{\text{pro\'et}}(X_0,\xi)$ is not continuous with respect to the discrete topology on $S$. Thus, we are not able to generalize the construction of $Y_0$ to any semi-stable curve $X_0$. We overcome this issue by associating with each continuous representation of $\pi_1^{\text{pro\'et}}(X_0,\xi)$ a specific geometric covering of $X_0$.

\begin{defn}\label{formalcovering}Let $\widehat{X}$ be the completion of $X$ along $X_0$, then we denote by $\text{Et}_{\widehat{X}}$ the category of formal schemes that are \'etale over $\widehat{X}$. We define $\text{Cov}_{\widehat{X}}$ to be the full subcategory $\text{Et}_{\widehat{X}}$ given by the essential image of $\text{Cov}_{X_0}$ via the equivalence in \cite[Exp. IX Prop 1.7]{SGA1}. We call the objects of $\text{Cov}_{\widehat{X}}$ \textit{geometric coverings of $\widehat{X}$}.
\end{defn}

\begin{rmk}Note that the categories $\text{Cov}_{X_0}$ and $\text{Cov}_{X}$ are, in general, not equivalent.
A counterexample is given by stable curves over $S$ with smooth generic fibre and degenerate closed fibre. If $X$ is such a curve, then, by \cite[Prop.10.3.15]{Liu}, $X$ is a normal scheme and, by \hyperref[geomuni]{Proposition \ref*{geomuni}}, $\pi^{\text{pro\'et}}_1(X)$ is a profinite group. While, by  \hyperref[genlemmapro]{Proposition \ref*{genlemmapro}}, $\pi^{\text{pro\'et}}_1(X_0)\simeq\mathbb{Z}^{\star r}$ with $r=p_a(X_0)\geq 2$, so the groups $\pi^{\text{pro\'et}}_1(X)$ and $\pi^{\text{pro\'et}}_1(X_0)$ are not isomorphic.

This counterexample also shows that there isn't a specialization morphism between the topological groups 
$\pi_1^{\text{pro\'et}}(X_{\widebar{K}})$ and $\pi_1^{\text{pro\'et}}(X_0)$  that lifts the \'etale specialization map. Indeed, any continuous morphism 
\[\text{sp}\colon\pi_1^{\text{pro\'et}}(X_{\widebar{K}})\simeq\pi_1^{\text{\'et}}(X_{\widebar{K}})\to\ZZ^{\star r}\]
factors through a finite quotient of $\pi_1^{\text{\'et}}(X_{\widebar{K}})$. Since $\ZZ^{\star r}$ is a free group, this implies that $\text{sp}$ is the zero map. Under the same assumptions, the \'etale specialization map is surjective, hence these morphisms are not compatible.
\end{rmk}

\begin{lem}\label{repcts}Let $\xi$ be a geometric point of $X_0$ and, for $j=1,\dots,N$, let $C_j$ be the irreducible components of $X_0$ and $\widebar{C_j}$ their normalizations, then there exists an equivalence of categories  
\[\Rp^{\text{cts}}_K(\pi^{\text{pro\'et}}_1(X_0,\xi)) \simeq\Rp^{\text{cts}}_K(\mathbb{Z}^{\star |I|-N+1 }\star\pi^{\textnormal{\'et}}_1(\widebar{C_1})\star\dots\star \pi^{\textnormal{\'et}}_1(\widebar{C_N})).\] 
\end{lem}
\proof Setting $r=|I|-N+1$, we consider the fixed isomorphism, 
\[\alpha\colon\pi^{\text{pro\'et}}_1(X_0,\xi)\simeq\mathbb{Z}^{\star r }\star_N\pi^{\text{\'et}}_1(\widebar{C_1})\star_N\dots\star_N \pi^{\text{\'et}}_1(\widebar{C_N}),\]
whose existence was proved in \hyperref[genlemmapro]{Proposition \ref*{genlemmapro}}.

By \hyperref[noohicoproductcompl]{Corollary \ref*{noohicoproductcompl}},  $\mathbb{Z}^{\star r}\star_N\pi^{\text{\'et}}_1(\widebar{C_1})\star_N\dots\star_N \pi^{\text{\'et}}_1(\widebar{C_N})$ is isomorphic to the Ra\u{\i}kov completion of  $\mathbb{Z}^{\star r}\star_{\mathcal{B}}\pi^{\text{\'et}}_1(\widebar{C_1})\star_{\mathcal{B}}\dots\star_{\mathcal{B}} \pi^{\text{\'et}}_1(\widebar{C_N})$, defined as in  \hyperref[basisneigh]{Lemma \ref*{basisneigh}}.
To simplify the notation let $\pi_{\mathcal{B}}=\mathbb{Z}^{\star r}\star_{\mathcal{B}}\pi^{\text{\'et}}_1(\widebar{C_1})\star_{\mathcal{B}}\dots\star_{\mathcal{B}} \pi^{\text{\'et}}_1(\widebar{C_N})$.
By \cite[Thm. 3.6.10]{AT}, there exists a continuous morphism $\sigma\colon \pi_{\mathcal{B}}\to\pi^{\text{pro\'et}}_1(X_0,\xi)$, whose image is dense. Hence $\sigma$ induces a fully faithful functor
\[\widetilde{\sigma}\colon\Rp^{\text{cts}}_K(\pi^{\text{pro\'et}}_1(X_0,\xi)) \to\Rp^{\text{cts}}_K(\pi_{\mathcal{B}}).\]
Let $(V,\rho)$ be a continuous representation of  $\pi_{\mathcal{B}}$, then, by \hyperref[repdiscrete]{Lemma \ref*{repdiscrete}}, $\rho$ induces a morphism $\rho\colon \pi_{\mathcal{B}}\to\text{Aut}(V)$ that is continuous with respect to the discrete topology on $\text{Aut}(V)$. Since groups with discrete topology are Ra\u{\i}kov complete, by \cite[Prop. 3.6.12]{AT}, $\rho$ admits an extension to $\hat{\rho}\colon \pi^{\text{pro\'et}}_1(X_0,\xi)\to\text{Aut}(V)$ such that $ \hat{\rho}\circ \sigma=\rho$. This implies that $\widetilde{\sigma}$ is an equivalence of categories.

Futhermore, as in \hyperref[noohicoproductcompl]{Corollary \ref*{noohicoproductcompl}}, we see that the identity map induces an equivalence of categories
\[\Rp^{\text{cts}}_K(\pi_{\mathcal{B}}) \simeq\Rp^{\text{cts}}_K(\mathbb{Z}^{\star r}\star\pi^{\text{\'et}}_1(\widebar{C_1})\star\dots\star \pi^{\text{\'et}}_1(\widebar{C_N})).\] 
Composing this functor with $\widetilde{\sigma}$, we construct the desired equivalence of categories.
\endproof

For the remaining of this article, we fix an equivalence of categories, as constructed in the above lemma.

Let us consider an element $(V,\rho)\in \Rp^{\text{cts}}_K(\pi^{\text{proet}}_1(X_0,\xi))$, then $(V,\rho)$ corresponds, via the equivalence of categories constructed in the proof of \hyperref[repcts]{Lemma \ref*{repcts}}, to a $K$-linear representation 
\[\rho\colon\mathbb{Z}^{\star r }\star\pi^{\text{\'et}}_1(\widebar{C_1})\star\dots\star \pi^{\text{\'et}}_1(\widebar{C_N})\to\text{Aut}(V),\]
which, by  \hyperref[repdiscrete]{Lemma \ref*{repdiscrete}}, is continuous with respect to the discrete topology on $\text{Aut}(V)$. 
Thus, by the universal property of the free product, $(V,\rho)$ corresponds to the following data:
 \begin{itemize}
 \item a continuous morphism $\rho^{\text{dis}}_i\colon\mathbb{Z}\to\text{Aut}(V)$ for $i=1,\dots,r$,
 \item a continuous morphism $\rho^{\text{\'et}}_j\colon\pi^{\text{\'et}}_1(\widebar{C_j})\to\text{Aut}(V)$ for $j=1,\dots,N$.
 \end{itemize} 
By continuity, each morphism   $\rho^{\text{\'et}}_j$ factors through a finite quotient of $\pi^{\text{\'et}}_1(\widebar{C_j})$, which we call $G_i$.
In particular, by the universal property of the free product, $(V,\rho)$  factors through a continuous $K$-linear representation
 \[\rho\colon\mathbb{Z}^{\star r}\star G_1\star\cdots\star G_N\to  \text{Aut}(V).\]

Clearly, $\mathbb{Z}^{\star r}\star G_1\star\cdots\star G_N$ is  a quotient of $\mathbb{Z}^{\star r }\star\pi^{\text{\'et}}_1(\widebar{C_1})\star\dots\star \pi^{\text{\'et}}_1(\widebar{C_N})$. Since it is a discrete group, by \cite[Prop. 3.6.12]{AT} it is also a quotient of  $\pi^{\text{pro\'et}}_1(X_0,\xi)$ and we denote the quotient map by $q$.

By \hyperref[proet]{Proposition \ref*{proet}}, the set $\mathbb{Z}^{\star r}\star G_1\star\cdots\star G_N$, endowed with the action given by $q$,  corresponds to a connected geometric covering  of $X_0$, which we denote by $Y_0^{\rho}$.

\begin{defn}\label{yrho} We set ${\mathcal{Y}}_{\rho}$ to be the  geometric covering of $\widehat{X}$ that corresponds to the geometric covering $Y_0^{\rho}$ of $X_0$ defined above.
\end{defn}  

By construction, we see that
\begin{equation}\label{autY}
\text{Aut}({\mathcal{Y}}_{\rho}|\widehat{X})\simeq\text{Aut}(Y_0^{\rho}|X_0)\simeq(\mathbb{Z}^{\star r}\star G_1\star\cdots\star G_N)^{\text{op}}.
\end{equation}

Similarly, we can endow the set $G_1\times\cdots\times G_N$ with a $\pi^{\text{pro\'et}}_1(X_0)$-action, by composing the map $q$ with the quotient map 
\[\alpha\colon\mathbb{Z}^{\star r}\star G_1\star\cdots\star G_N\to G_1\times\cdots\times G_N.\]
Hence, we can associate with $ G_1\times\cdots\times G_N$ a finite \'etale cover $Z_0^{\rho}$  of $X_0$.

\begin{defn}\label{zrho} We set ${\mathcal{Z}}_{\rho}$ to be the finite \'etale covering of $\widehat{X}$ that corresponds to the finite \'etale covering $Z_0^{\rho}$ of $X_0$ defined above.
\end{defn} 

We can observe  that
\begin{equation}\label{autZ}
\text{Aut}(\mathcal{Z}_{\rho}|\widehat{X})\simeq\text{Aut}(Z_0^{\rho}|X_0)\simeq (G_1\times\cdots\times G_N)^{\text{op}}.
\end{equation}

Moreover, $\mathcal{Y}_{\rho}\to\widehat{X}$ factors through $q\colon\mathcal{Y}_{\rho}\to\mathcal{Z}_{\rho}$ and we  have 
\begin{equation}\label{auta}
\text{Aut}(\mathcal{Y}_{\rho}|\mathcal{Z}_{\rho})\simeq\text{Aut}(Y_0^{\rho}|Z_0^{\rho})\simeq\text{ker}(\alpha)^{\text{op}}.
\end{equation}

Note that the morphism $\mathcal{Y}_{\rho}\to \widehat{X}\to X\to S$ corresponds to  $A\to \Gamma(\mathcal{Y}_{\rho},\OO_{\mathcal{Y}_{\rho}})$. Hence, a coherent $\OO_{\mathcal{Y}_{\rho}}$-module is a sheaf of $A$-modules.
 
\begin{defn}\label{defmerom} Given  $\FF$ a coherent sheaf on  $\mathcal{Y}_{\rho}$,  we call \textit{meromorphic descent data relative to $\mathcal{Z}_{\rho}$ on $\FF$} a collection of elements 
\[h_{w}\in\hm^0(\mathcal{Y}_{\rho}, \textit{Hom}_{\mathcal{O}_{\mathcal{Y}_{\rho}}}(\FF,w^{\star}\FF)\otimes_A K),\text{ }w\in \text{Aut}(\mathcal{Y}_{\rho}|\mathcal{Z}_{\rho})\]
that satisfy:
\begin{itemize}
\item  the co-cycle condition: $w^{\star}h_{w'}\circ h_w=h_{w'\circ w}\text{ for every }w,w'\in \text{Aut}(\mathcal{Y}_{\rho}|\mathcal{Z}_{\rho})$;
\item the identity condition: $h_{\text{Id}}=\text{Id}_{\FF\otimes_A K}$.
\end{itemize}
\end{defn}

\begin{defn}\label{morphmerom} Given $\{\FF,h_w\}_{w\in\text{Aut}(\mathcal{Y}_{\rho}|\mathcal{Z}_{\rho})}$ and $\{\GG,k_w\}_{w\in\text{Aut}(\mathcal{Y}_{\rho}|\mathcal{Z}_{\rho})}$ two coherent sheaves on $\mathcal{Y}_{\rho}$ with meromorphic descent data  relative to $\mathcal{Z}_{\rho}$,  a \textit{morphism of meromorphic  descent data} from $\{\FF,h_w\}$ to $\{\GG,k_w\}$ is given by an element 
\[f\in\hm^0(\mathcal{Y}_{\rho},\textit{Hom}_{\mathcal{O}_{\mathcal{Y}_{\rho}}}(\FF,\GG)\otimes_A K)\] such that for every $w\in \text{Aut}(\mathcal{Y}_{\rho}|\mathcal{Z}_{\rho})$
\[k_w\circ f=w^{\star}(f)\circ h_w.\]

We denote by $\text{Coh}^m(\mathcal{Y}_{\rho}|\mathcal{Z}_{\rho})$ the category of coherent sheaves on $\mathcal{Y}_{\rho}$ with meromorphic descent data relative to $\mathcal{Z}_{\rho}$.
\end{defn}

\begin{defn} \label{defdescent}Let $\{\FF,h_w\}_{w\in\text{Aut}(\mathcal{Y}_{\rho}|\mathcal{Z}_{\rho})}$ be a coherent sheaf on ${\mathcal{Y}_{\rho}}$ with meromorphic descent data relative to ${\mathcal{Z}_{\rho}}$, we say that $\{\FF,h_w\}$ \textit{descends to a coherent sheaf on ${\mathcal{Z}_{\rho}}$} if there exists $\GG\in\text{Coh}({\mathcal{Z}_{\rho}})$ such that 
\[\{\FF,h_w\}_{w\in\text{Aut}(\mathcal{Y}_{\rho}|\mathcal{Z}_{\rho})}\simeq\{q^{\star}\GG,h^q_w\}_{w\in\text{Aut}(\mathcal{Y}_{\rho}|\mathcal{Z}_{\rho})},\] 
where $h^q_w\colon q^{\star}\GG\to w^{\star}q^{\star}\GG$ are the natural isomorphisms.
\end{defn}

The following proposition is a generalization of \cite[Lemma 1]{Gies1}.

\begin{prop} \label{infdesc1}For every coherent sheaf $\{\FF,h_w\}_{w\in\text{Aut}(\mathcal{Y}_{\rho}|\mathcal{Z}_{\rho})}$ with meromorphic descent  data  relative to $\mathcal{Z}_{\rho}$, there exists a coherent sheaf $\{\FF',k_w\}_{\Aut(\mathcal{Y}_{\rho}|\mathcal{Z}_{\rho})}$   with meromorphic descent data relative to $\mathcal{Z}_{\rho}$ that is isomorphic to $\{\FF,h_w\}_{w\in\text{Aut}(\mathcal{Y}_{\rho}|\mathcal{Z}_{\rho})}$ and such that
\[k_w\in\text{H}^0(X,\text{Hom}_{\OO_{{\mathcal{Y}}_{\rho}}}(\FF',w^{\star}\FF')).\]
\end{prop}
\proof 
As in \cite[Lemma 1]{Gies1}, it suffices to show that, for any $\text{Aut}({\mathcal{Y}_{\rho}}|{\mathcal{Z}_{\rho}})$-invariant open  $U\subsetneq \mathcal{Y}_{\rho}$, there exists a quasi-compact open $V$ of ${\mathcal{Y}}_{\rho}$ such that
\begin{itemize}
\item $V$ is not contained in $U$,
\item $V\cap wV\subseteq U$ for all $w\in \text{Aut}(\mathcal{Y}_{\rho}|\mathcal{Z}_{\rho})$, $w\neq Id_{\mathcal{Y}_{\rho}}$.
\end{itemize}

Let $\xi_j\in C_j$ be a fixed geometric point for every $j=1,\dots,N$. By \hyperref[genlemmapro]{Proposition \ref*{genlemmapro}} we can choose an irreducible component $\mathcal{Y}^j_{\emptyset}$ of $\mathcal{Y}_{\rho}$, such $F_{\xi_j}(\mathcal{Y}^j_{\emptyset})=G_j$.  
Given a word $s\in\mathbb{Z}^{\star r}\star G_1\star\cdots\star G_N$, we denote $\mathcal{Y}^j_s$ the irreducible component $\mathcal{Y}^j_s:=s(\mathcal{Y}^j_{\emptyset})$, which corresponds via the functor $F_{\xi_j}$ to the $G_{j}$-orbit of $s$. By \hyperref[genlemmapro]{Proposition \ref*{genlemmapro}}, the set $\{\mathcal{Y}^j_s\}_{s,j}$ contains the set of all irreducible components of $\mathcal{Y}_{\rho}$.
Since the action of $\text{ker}(\alpha)^{\text{op}}$ on $\mathbb{Z}^{\star r}\star G_1\star\cdots\star G_N$ is defined by right concatenation, given $\mathcal{Y}^j_s$ an irreducible component of $\mathcal{Y}_{\rho}$ and $w\in\text{ker}(\alpha)^{\text{op}}$, $w\neq Id_{\mathcal{Y}_{\rho}}$, we have
\[w(G_js)=G_jsw\neq G_js\text{ and }w(\mathcal{Y}^j_s)=\mathcal{Y}^j_{sw}\neq \mathcal{Y}^j_s.\] 
Hence, the action of $\text{ker}(\alpha)^{\text{op}}$ on the set of irreducible components of $\mathcal{Y}_{\rho}$ is free.

Let us suppose  that we are given an open $\text{Aut}({Y_{\rho}}|{Z_{\rho}})$-invariant set $U\subset \mathcal{Y}_{\rho}$, then for the construction of $V$ there are two possible cases.

\noindent\textit{First case}: there exists $x\in {\mathcal{Y}}_{\rho}\setminus U$ that is  a non-singular point.

Then we set $\mathcal{Y}_x$ to be the irreducible component of ${\mathcal{Y}}_{\rho}$ containing $x$,  $I_x$ to be the set of the singular points of $\mathcal{Y}_x$ and we define $V= \mathcal{Y}_x\setminus I_x$.
By construction, $V$ is not contained in $U$ and we have, for all $w\in\text{ker}(\alpha)^{\text{op}}$, $w\neq Id_{\mathcal{Y}_{\rho}}$,
\[V\cap wV=\emptyset\subset U.\]

\noindent\textit{Second case}: $\mathcal{Y}_{\rho}\setminus U\subset I$, where $I$ is the set of singular points of $\mathcal{Y}_{\rho}$.

Let $x\in \mathcal{Y}_{\rho}\setminus U$, then $x$ belongs to exactly two irreducible components of ${\mathcal{Y}_{\rho}}$, say $\mathcal{Y}^i_s$ and $\mathcal{Y}^l_t$. Let $I_x$ be the set of singular points of $\mathcal{Y}^i_s\cup Y^l_t$ different from $x$, then we set  $V=(\mathcal{Y}^i_s\cup \mathcal{Y}^l_t)\setminus I_x$.
Clearly, $V$ is not contained in $U$. Moreover, if $w\in\text{ker}(\alpha)^{\text{op}}$ is a non-trivial word, then 
\[V\cap wV=((\mathcal{Y}^i_s\cap \mathcal{Y}^l_{tw})\cup(\mathcal{Y}^l_t\cap \mathcal{Y}^i_{sw}))\setminus\{\text{sing. pts}\}.\] 
Thus, there are three possibilities:
\begin{itemize}
\item $\mathcal{Y}^l_{tw}\neq \mathcal{Y}^i_s$ and $\mathcal{Y}^i_{sw}\neq \mathcal{Y}^l_t $, that implies $V\cap wV=\emptyset\subset U$,
\item $\mathcal{Y}^l_{tw}=\mathcal{Y}^i_s$ and $\mathcal{Y}^i_{sw}\neq \mathcal{Y}^l_t $, that implies $V\cap wV=\mathcal{Y}^i_s\setminus \{\text{sing. pts of } \mathcal{Y}^i_s\}\subset U$,
\item $\mathcal{Y}^l_{tw}\neq \mathcal{Y}^i_s$ and $\mathcal{Y}^i_{sw}=\mathcal{Y}^l_t $, that implies $V\cap wV=\mathcal{Y}^l_t\setminus \{\text{sing. pts of } \mathcal{Y}^l_t\}\subset U$.
\end{itemize}
Note that the case where $\mathcal{Y}^l_{tw}=\mathcal{Y}^i_s$ and $\mathcal{Y}^i_{sw}=\mathcal{Y}^l_t $ does not occur because it would imply that $w^2=Id_{\mathcal{Y}_{\rho}}$, which is not possible because $\text{ker}(\alpha)^{\text{op}}$ is torsion free. 
\endproof

\begin{rmk}\label{nodirectdescent} The action of $(\mathbb{Z}^{\star r}\star G_1\star\cdots\star G_N)^{\text{op}}$ on the set of irreducible components of $\mathcal{Y}_{\rho}$ is not free. Indeed, if $\emptyset$ is the empty word  and  $\mathcal{Y}^j_{\emptyset}$  is the irreducible component of $\mathcal{Y}_{\rho}$ that corresponds to $G_j\subset F_{\xi}(\mathcal{Y}_{\rho})$, then  for every $g_j\in G_J$,
\[g_j(\mathcal{Y}^j_{\emptyset})=\mathcal{Y}^j_{g_j}=\mathcal{Y}^j_{\emptyset}.\]
\end{rmk}

The following theorem generalizes \cite[Lemma 2]{Gies1}.

\begin{thm} \label{infdesc}Any coherent sheaf $\{\FF,h_w\}_{w\in\text{Aut}({\mathcal{Y}_{\rho}}|{\mathcal{Z}_{\rho}})}$ on $\mathcal{Y}_{\rho}$ with meromorphic descent data  relative to $\mathcal{Z}_{\rho}$ descends to a coherent sheaf on $\mathcal{Z}_{\rho}$.
\end{thm}
\proof As in \cite[Lemma 2]{Gies1},  it suffices to prove that there exists a quasi-compact open subscheme $T$ of $\mathcal{Y}_{\rho}$ such that its $\text{Aut}({\mathcal{Y}_{\rho}}|{\mathcal{Z}_{\rho}})$-translates  cover ${\mathcal{Y}}_{\rho}$. 

We fix  a non-trivial word $w\in\text{ker}(\alpha)^{\text{op}}$. 
Note that the irreducible components of the form $\mathcal{Y}^j_{s'}$ with $\alpha(s')=\alpha(s)$, defined as in the previous theorem's proof, are $\text{ker}(\alpha)^{\text{op}}$-translates of $\mathcal{Y}^j_{ws}$. Indeed, the word 
$t=s^{-1}w^{-1}s'$ satifies 
\[t(\mathcal{Y}^j_{ws})=\mathcal{Y}^j_{wst}=\mathcal{Y}^j_{s'}.\] 

Given an element $g=(g_1,\dots,g_N)\in G_1\times\cdots\times G_N$, we denote by $\sigma(g)$ the word $g_1\cdots g_N\in \mathbb{Z}^{\star r}\star G_1\star\cdots\star G_N$ with letters $g_i\in G_i$ and we define the map
\[\sigma\colon G_1\times\cdots\times G_N\to \mathbb{Z}^{\star r}\star G_1\star\cdots\star G_N\text{, }\sigma(g_1,\dots,g_N)=g_1\cdots g_N.\] 
We denote by $1_i$ the word whose only letter is the element $1\in\mathbb{Z}$ belonging to the $i$-th copy of $\ZZ$. Then we set 
\[T_G=\bigcup_{j=1}^N\bigcup_{g\in G_1\times\cdots\times G_N} \bigg(\mathcal{Y}^j_{w\sigma(g)}\cup\bigcup_{i=1}^r \mathcal{Y}^j_{1_iw\sigma(g)}\bigg).\]
and we define $I_G$ to be set of points of $T_G$ that are intersection points with  irreducible components of ${\mathcal{Y}_{\rho}}$ not contained in $T_G$. For every $x\in I_G$, let $V_x$ be a quasi-compact open neighborhood of $x$, then we set 
\[T=T_G\setminus I_G\cup \bigcup_{x\in I_G} V_x.\] 

Since $I_G$ is a finite set, $T$ is by construction a quasi-compact open of $\mathcal{Y}_{\rho}$. Thus, it suffices to prove that its $\text{ker}(\alpha)^{\text{op}}$-translates cover  ${\mathcal{Y}}_{\rho}$.

Given $s\in  \mathbb{Z}^{\star r}\star G_1\star\cdots\star G_N$, we set $g_s:=\alpha(s)$.
Since $\alpha(s)=\alpha(w\sigma(g_s))$, there exists $t\in\text{ker}(\alpha)^{\text{op}}$ such that 
\[t(\mathcal{Y}^j_{w\sigma(g_s)})=\mathcal{Y}^j_{s}.\]
This implies that
\[{\mathcal{Y}_{\rho}}=\bigcup_{t\in \text{ker}(\alpha)^{\text{op}}} t(T).\]

\endproof

\begin{thm}\label{equivalence}  Given $(V,\rho)\in \Rp^{\textnormal{cts}}_K(\pi^{\text{pro\'et}}_1(X_0,\xi))$, let $Z_{\rho}$ be the finite \'etale covering of $X$ corresponding to $\mathcal{Z}_{\rho}$  and $Z_K^{\rho}$ its generic fibre, then  the category  $\text{Coh}^m(\mathcal{Y}_{\rho}|\mathcal{Z}_{\rho})$ of coherent sheaves on $\mathcal{Y}_{\rho}$ with meromorphic descent relative to $\mathcal{Z}_{\rho}$  is equivalent to the category $\text{Coh}(Z_K^{\rho})$ of coherent sheaves on $Z_K^{\rho}$.
\end{thm}
\proof By \hyperref[infdesc]{Theorem \ref*{infdesc}} and \cite[Prop. 1]{Gies1},   $\text{Coh}^m(\mathcal{Y}_{\rho}|\mathcal{Z}_{\rho})$  is equivalent to the category $\text{Coh}^K(\mathcal{Z}_{\rho})$, whose objects are coherent sheaves on $\mathcal{Z}_{\rho}$ and whose morphisms defined by 
\[\text{Hom}_{\text{Coh}^K(\mathcal{Z}_{\rho})}(\FF,\GG):=\text{Hom}_{\mathcal{O}_{\mathcal{Z}_{\rho}}}(\FF,\GG)\otimes_A K.\] 

Moreover, by Grothendieck's existence theorem \cite[Cor.5.1.6]{EGAIII}, the category $\text{Coh}^K(\mathcal{Z}_{\rho})$ is equivalent to the category $\text{Coh}^K( {Z}_{\rho})$, whose objects are coherent sheaves on $ {Z}_{\rho}$ and whose maps are  given by
\[\text{Hom}_{\text{Coh}^K( {Z}_{\rho})}(\FF,\GG):=\text{Hom}_{\mathcal{O}_{{Z}_{\rho}}}(\FF,\GG)\otimes_A K.\] 

Denoting $j\colon Z_K^{\rho}\to Z_{\rho}$ the open immersion, it suffices to show that the functor 
\[j^{\star}\colon\text{Coh}^K(Z_{\rho})\to\text{Coh}(Z_K^{\rho})\]
is an equivalence of categories.
By flat base change \cite[5.2.27]{Liu}, for every coherent sheaf $\FF$ on $Z_{\rho}$ and for any $p\geq 0$,
\[\text{H}^p(Z_{\rho},\FF)\otimes_A K\cong \text{H}^p(Z_K^{\rho},j^{\star}\FF).\]
Applying this  to the sheaf $\textit{Hom}_{\mathcal{O}_{Z_{\rho}}}(\FF,\GG)$, for every $\FF$ and $\GG$ coherent sheaves, we get that $j^{\star}$ is a fully faithful functor.
Moreover, since $Z_{\rho}$ is proper over $S$, by \cite[Thm. 9.4.8]{EGAI} the functor $j^{\star}$ is essentially surjective.
\endproof

\begin{rmk}\label{morphK} Let $\text{Coh}^K(\widehat{X})$ be the category, whose objects are coherent sheaves on $\widehat{X}$ and whose morphisms defined by 
\[\text{Hom}_{\text{Coh}^K(\widehat{X})}(\FF,\GG):=\text{Hom}_{\mathcal{O}_{\widehat{X}}}(\FF,\GG)\otimes_A K.\] 
From the same reasoning of the previous result's proof, it follows that the category $\text{Coh}^K(\widehat{X})$ 
is equivalent to the category of coherent sheaf on $X_K$. 
\end{rmk}

We prove now that meromorphic data for a coherent sheaf $\FF$ on  $\mathcal{Y}_{\rho}$ descend to a coherent sheaf on $X_K$.

\begin{defn} Extending \hyperref[defmerom]{Definition \ref*{defmerom}}, we define \textit{meromorphic descent data for a coherent sheaf $\FF$} on  $\mathcal{Y}_{\rho}$,  to be a collection of elements 
\[h_w\in\hm^0(\mathcal{Y}_{\rho}, \textit{Hom}_{\mathcal{O}_{\mathcal{Y}_{\rho}}}(\FF,w^{\star}\FF)\otimes_A K),\text{ }w\in \text{Aut}(\mathcal{Y}_{\rho}|\widehat{X})\]
that satisfy:
\begin{itemize}
\item  the co-cycle condition: $w^{\star}h_{w'}\circ h_w=h_{w'\circ w}\text{ for every }w,w'\in \text{Aut}(\mathcal{Y}_{\rho}|\widehat{X})$;
\item  the identity condition: $h_{\text{Id}}=\text{Id}_{\FF\otimes_A K}$.
\end{itemize}

The definition of morphisms between coherent sheaves on $\mathcal{Y}_{\rho}$ with meromorphic descent data is analogous to \hyperref[morphmerom]{Definition \ref*{morphmerom}}.

We denote by $\text{Coh}^m(\mathcal{Y}_{\rho}|\widehat{X})$ the category of coherent sheaves on $\mathcal{Y}_{\rho}$ with meromorphic descent data.
\end{defn}

\begin{rmk} Let $\FF$ be a coherent sheaf on $\mathcal{Y}_{\rho}$, then the meromorphic descent data $(\FF,h_w)_{w\in \text{Aut}(\mathcal{Y}_{\rho}|\widehat{X})}$ induces in particular meromorphic descent data $(\FF,h_w)_{w\in \text{Aut}(\mathcal{Y}_{\rho}|\mathcal{Z}_{\rho})}$ relative to $\mathcal{Z}_{\rho}$. 

By \hyperref[equivalence]{Theorem \ref*{equivalence}}, this implies that $(\FF,h_w)_{w\in \text{Aut}(\mathcal{Y}_{\rho}|\mathcal{Z}_{\rho})}$ descends to a coherent sheaf on $\mathcal{Z}_{\rho}$.
\end{rmk}

\begin{lem}\label{gF} Let $(\FF,h_w)_{w\in \text{Aut}(\mathcal{Y}_{\rho}|\widehat{X})}$ be a coherent sheaf with meromorphic descent data, which by \hyperref[equivalence]{Theorem \ref*{equivalence}} descends to a coherent sheaf $\FF_Z$ on $Z_K^{\rho}$, then for every $g\in \mathbb{Z}^{\star r}\star G_1\star\cdots\star G_N$ the coherent sheaf $(g^{\star}\FF,g^{\star}h_{g\circ w \circ g^{-1}})_{w\in \text{Aut}(\mathcal{Y}_{\rho}|\mathcal{Z}_{\rho})}$ descends to ${\alpha(g)}^{\star}\FF_Z$ on $\mathcal{Z}_{\rho}$.
\end{lem}
\proof
Let $q_Z\colon \mathcal{Y}_{\rho} \to \mathcal{Z}_{\rho}$ be the constructed geometric covering. Since the sheaf $(\FF,h_w)_{w\in \text{Aut}(\mathcal{Y}_{\rho}|\mathcal{Z}_{\rho})}$  descends to a coherent sheaf $\FF_Z$  on $\mathcal{Z}_{\rho}$, there exists an isomorphism $\psi$ 
such that, for every $w'\in\text{ker}(\alpha)^{\text{op}}$, the following diagram commutes
\begin{equation}\label{fhw}
\begin{CD}
q_Z^{\star}\FF_Z\otimes_A K@>\psi>> \FF\otimes_A K\\
@V{id}VV @VV{ h_{w'}}V\\
{w'}^{\star}q_Z^{\star}\FF_Z\otimes_A K=q_Z^{\star}\FF_Z\otimes_A K@>w'^{\star}\psi>> {w'}^{\star}\FF\otimes_A K\text{ .}\\
\end{CD}
\end{equation}

Moreover, for every $g\in \mathbb{Z}^{\star r}\star G_1\star\cdots\star G_N$, $\alpha(g)\circ q_Z=q_Z\circ g$ and thus
\[
g^{\star}q_Z^{\star}\FF_Z=q_Z^{\star}{\alpha(g)}^{\star}\FF_Z\otimes_A K. 
\]

Finally, if we take $w'=g\circ w\circ g^{-1}$ and we apply  $g^{\star}$ to \eqref{fhw}, we  see that the following diagram commutes
\begin{equation*}\begin{CD}
q_Z^{\star}{\alpha(g)}^{\star}\FF_Z\otimes_A K@>g^{\star}\psi>> g^{\star}\FF\otimes_A K\\
@V{id}VV @VV{g^{\star}h_{g\circ w\circ g^{-1}}}V \\
w^{\star}q_Z^{\star}{\alpha(g)}^{\star}\FF_Z\otimes_A K@>w^{\star}g^{\star}\psi>> w^{\star}g^{\star}\FF\otimes_A K\text{ .}\\
\end{CD}
\end{equation*}
Therefore, we can conclude that  $(g^{\star}\FF,g^{\star}h_{g\circ w \circ g^{-1}})_{w\in \text{Aut}(\mathcal{Y}_{\rho}|\mathcal{Z}_{\rho})}$ descends to the coherent sheaf ${\alpha(g)}^{\star}\FF_Z$.
\endproof

\begin{thm} \label{findesc} Let $q_X\colon \mathcal{Y}_{\rho}\to\widehat{X}$ be the constructed geometric covering, then the pullback $q_X^{\star}$ induces an equivalence of categories between the category $\text{Coh}^K(\widehat{X})$, defined as in \hyperref[morphK]{Remark \ref*{morphK}}, and the category $\text{Coh}^m(\mathcal{Y}_{\rho}|\widehat{X})$ of coherent sheaves on $\mathcal{Y}_{\rho}$ with meromorphic descent data. 

In particular, the category $\text{Coh}^m(\mathcal{Y}_{\rho}|\widehat{X})$  is equivalent to the category $\text{Coh}(X_K)$ of coherent sheaves on $X_K$. 
\end{thm}
\proof 

Clearly, the pullback of a coherent sheaf on $\widehat{X}$ along $q_X\colon \mathcal{Y}_{\rho}\to \widehat{X}$ can be endowed with natural meromorphic descent data. This construction induces the desired functor
\[q_X^{\star}\colon \text{Coh}^K(\widehat{X}) \to \text{Coh}^m(\mathcal{Y}_{\rho}|\widehat{X}).\]

Let $\text{Coh}^m(\mathcal{Z}_{\rho}|\widehat{X})$ be the category of coherent sheaves $\FF$ on $\mathcal{Z}_{\rho}$ with meromorphic descent data $\{h_g\}_{g\in\text{Aut}(\mathcal{Z}_{\rho}|\widehat{X})}$. By construction, the functor $q_X^{\star}$ factors as follows
\[
\begin{tikzcd}
 & \text{Coh}^m(\mathcal{Z}_{\rho}|\widehat{X}) \arrow{dr}{q_Z^{\star}} \\
\text{Coh}^K(\widehat{X})\arrow{ur}{q_{Z|X}^{\star}} \arrow{rr}{q_X^{\star}} && \text{Coh}^m(\mathcal{Y}_{\rho}|\widehat{X}) ,
\end{tikzcd}
\]  
where $q_{Z|X}\colon \mathcal{Z}_{\rho}\to \widehat{X}$ and $q_Z\colon \mathcal{Y}_{\rho} \to \mathcal{Z}_{\rho}$ are the geometric coverings we defined.

We recall that, as explained in \hyperref[morphK]{Remark \ref*{morphK}}, the category $\text{Coh}^K(\widehat{X})$ is naturally equivalent to the category of coherent sheaves on $X_K$. 
Furthermore, the argument in the proof of \hyperref[equivalence]{Theorem \ref*{equivalence}} implies that the category  $\text{Coh}^m(\mathcal{Z}_{\rho}|\widehat{X})$ is equivalent to the category $\text{Coh}(Z_K^{\rho}|X_K)$ of coherent sheaves on $Z_K^{\rho}$ with descent data relative to the finite \'etale morphism $p\colon Z_K^{\rho}\to X_K$ and the following diagram commutes
\[
\begin{CD}
\text{Coh}^K(\widehat{X})@>>>\text{Coh}(X_K)\\
@V{q_{Z|X}^{\star}}VV @VV{p^{\star}}V\\
\text{Coh}^m(\mathcal{Z}_{\rho}|\widehat{X}) @>>>\text{Coh}(Z_K^{\rho}|X_K) \text{ }.
\end{CD}
\]
Since finite \'etale morphism are of effective descent, the functor $p^{\star}$ is an equivalence of categories, hence, so is $q_{Z|X}^{\star}$. To prove the theorem, it suffices to show that the functor $q_Z^{\star}$ is an equivalence of categories. 
We  prove first that $q_Z^{\star}$ is essentially surjective. 

Let $(\FF,h_w)_{w\in\text{Aut}(\mathcal{Y}_{\rho}|\widehat{X})}$ be a coherent sheaf with meromorphic descent data, then, by \hyperref[equivalence]{Theorem \ref*{equivalence}}, $(\FF,h_w)_{w\in \text{Aut}(\mathcal{Y}_{\rho}|\mathcal{Z}_{\rho})}$ descends to a coherent sheaf $\FF_Z$ on $\mathcal{Z}_{\rho}$. It remains to construct meromorphic descent data for the sheaf $\FF_Z$ relative to the map $q_{Z|X}\colon \mathcal{Z}_{\rho}\to \widehat{X}$. 
By \hyperref[gF]{Lemma \ref*{gF}}, the coherent sheaf with meromorphic descent data $(g^{\star}\FF,g^{\star}h_{g\circ w \circ g^{-1}})_{w\in \text{Aut}(\mathcal{Y}_{\rho}|\mathcal{Z}_{\rho})}$ descends to the sheaf ${\alpha(g)}^{\star}\FF_Z$ on $\mathcal{Z}_{\rho}$, for every $g\in \mathbb{Z}^{\star r}\star G_1\star\cdots\star G_N$. Hence, by \hyperref[equivalence]{Theorem \ref*{equivalence}}, we need to construct
\[h_{\alpha(g)}\in\text{Hom}(\FF_Z\otimes_A K,{\alpha(g)}^{\star}\FF_Z\otimes_A K)=\text{Hom}(\{\FF,h_w\},\{g^{\star}\FF,g^{\star}h_{g\circ w\circ g^{-1}}\}).\]
By the co-cycle condition of meromorphic descent data, the following diagram commutes
\[
\begin{CD}
\FF\otimes_A K@>{h_{g}}>> g^{\star}\FF\otimes_A K\\
@V{h_w}VV @VV{g^{\star}h_{g\circ w\circ g^{-1}}}V\\
\FF\otimes_A K@>>{w^{\star}h_{g}}>  g^{\star}\FF\otimes_A K\text{ .}\\
\end{CD}
\]
Hence, $h_g$ induces an isomorphism from $\FF_Z$ to $\alpha(g)^{\star}\FF_Z$, which only depends on $\alpha(g)$. Since $\{h_w\}_{w\in \text{Aut}(\mathcal{Y}_{\rho}|\widehat{X})}$ satisfy the co-cycle condition, so do the isomorphisms $\{h_{\alpha(g)}\}$. Therefore, the collection $\{h_{\alpha(g)}\}$ gives natural descent data on $\FF_Z$ relative to $q_{Z|X}\colon\mathcal{Z}_{\rho}\to\widehat{X}$.

By construction, there exists an isomorphism $\psi\colon q_Z^{\star}\FF_Z\otimes_A K \to \FF\otimes_A K $. Moreover, by construction of $\FF_Z$ and $h_{\alpha(w)}$, the following diagram commutes:
\[
\begin{CD}
q_Z^{\star}\FF_Z\otimes_A K@>{\psi}>> \FF\otimes_A K\\
@V{q_Z^{\star}h_{\alpha(w)}}VV @VV{h_{w}}V\\
q_Z^{\star}\FF_Z\otimes_A K @>>{w^{\star}\psi}>  w^{\star}\FF\otimes_A K\text{ .}\\
\end{CD}.
\]
Hence, $\psi$ is an isomorphism of coherent sheaves with meromorphic data and the functor $q_{Z}^{\star}$ is essentially surjective.

It remains to prove that the functor $q_{Z}^{\star}$ is fully faithful. Let $(\FF_Z,h_g)_{g\in \text{Aut}(\mathcal{Z}_{\rho}|\widehat{X})}$ and $(\GG_Z,k_g)_{g\in \text{Aut}(\mathcal{Z}_{\rho}|\widehat{X})}$ be coherent sheaves on $\mathcal{Z}_{\rho}$ with meromorphic descent data and let $(\FF,h_w)_{w\in \text{Aut}(\mathcal{Y}_{\rho}|\widehat{X})}$ and $(\GG,k_w)_{w\in \text{Aut}(\mathcal{Y}_{\rho}|\widehat{X})}$ be their pullback on $\mathcal{Y}_{\rho}$. 

Given two morphisms $f_1,f_2\colon (\FF_Z,h_g)\to (\GG_Z,k_g)$, if $q_Z^{\star}f_1=q_Z^{\star}f_2$ as morphisms of sheaves with meromorphic descent data, then they coincide in particular as morphisms of sheaves with meromorphic descent data relative to $\mathcal{Z}_{\rho}$. By \hyperref[equivalence]{Theorem \ref*{equivalence}}, this implies that $f_1=f_2$.

Let $f$ be a morphism between $(\FF,h_w)_{w\in \text{Aut}(\mathcal{Y}_{\rho}|\widehat{X})}$ and $(\GG,k_w)_{w\in \text{Aut}(\mathcal{Y}_{\rho}|\widehat{X})}$, then by \hyperref[equivalence]{Theorem \ref*{equivalence}} there exists a morphism of sheaves $\widebar{f}\colon(\FF_Z,h_g)\to (\GG_Z,k_g)$ such that $q_Z^{\star}\widebar{f}=f$. For every $g\in \text{Aut}(\mathcal{Z}_{\rho}|\widehat{X})$, there exists $s\in \text{Aut}(\mathcal{Y}_{\rho}|\widehat{X})$ such that $\alpha(s)=g$. The morphisms $k_g\circ \widebar{f}$ and $g^{\star}\widebar{f}\circ h_g$ correspond via \hyperref[equivalence]{Theorem \ref*{equivalence}} to $k_s\circ f$ and $s^{\star}f\circ h_s$, which coincide for every chosen $s$.  Hence, it is clear that $\widebar{f}$ is a morphism of meromorphic descent data.
\endproof

\section{Specialization functor}

In this section we construct the specialization functor between the category $\Rp^{\textnormal{cts}}_K(\pi^{\textnormal{pro\'et}}_1(X_0,\xi))$ and the category $\St(X_{\widebar K})$ of stratified bundles. We start recalling the definition and properties of the latter.

\begin{defn} Let $T$ be a smooth  scheme of finite type over a field $F$ of positive characteristic, $F_{T/F}^i$ the relative Frobenius and $T^{(i)}$ its $i$-th Frobenius twist, then an \textit{F-divided sheaf on $T$} is given by a sequence $(\EE_i,\sigma_i)_{i\geq 0}$,  where $\EE_i$ are bundles on $T^{(i)}$  and $\sigma_i \colon {F_{T/F}^i}^{\star}\EE_{i+1}\to\EE_i$ are $\mathcal{O}_{T^{(i)}}$-linear isomorphisms.
\end{defn}

\begin{defn}  Given $(\EE_i,\sigma_i)$ and $(\GG_i,\tau_i)$ F-divided sheaves on a scheme $T$ as above, a \textit{morphism of stratified bundles} from $(\EE_i,\sigma_i)$ to $(\GG_i,\tau_i)$ is defined as a sequence of $\mathcal{O}_{T^{(i)}}$-linear maps $\alpha=\{\alpha_i\colon\EE_i\to\GG_i\}$ such that the following diagram is commutative
\[
\begin{CD}
{F_{T/F}^i}^{\star}\EE_{i+1} @>{{F_{T/F}^i}^{\star}\alpha_{i+1}}>>{F_{T/F}^i}^{\star}\GG_{i+1}\\
@V\sigma_iVV @VV\tau_iV \\
\EE_i @>{\alpha_i}>> \GG_i\text{ .}\\
\end{CD}
\]
\end{defn}

\begin{defn} Let $T$ be a smooth  scheme of finite type over a field $F$ and $\DD_{T/F}$ the quasi coherent $\OO_T$-module of differential operators defined in \cite[Section 16]{EGA4}, then a \textit{stratified bundle on  $T$}  is a locally free $\OO_T$-module of finite rank endowed with a $\OO_T$-linear $\DD_{T/F}$-action extending the $\OO_T$-module structure via the inclusion $\OO_T\subset\DD_{T/F}.$ A \textit{morphism of stratified bundles} is a morphism of $\DD_{T/F}$-modules.
\end{defn}

\begin{thm}[\textbf{Katz's theorem}, \cite{Gies2}, Thm. 1.3] \label{Katzthm} Let $T$ be a smooth scheme of finite type over a perfect field $F$ of characteristic $p>0$, then the category of stratified bundles on $T$ and the category of F-divided sheaves on $T$ are equivalent. 
\end{thm}

If the base field is perfect, we will identify these two categories and we use the term stratified bundles for both definitions. Moreover, we will denote by $\St(T)$ the category of  stratified bundles on $T$.

\begin{prop}[\cite{Saa}, Section. VI.1.2]\label{saavedra} Let $T$ be a smooth scheme of finite type over a perfect field $F$, then the category $\St(T)$ of stratified bundles on $T$ is a rigid abelian tensor category.
Moreover, if $T$ has a rational point $x\in T(F)$,  the functor 
\[\omega_x\colon\St(T)\to\Vc_F \text{, }\omega_x(\EE_i,\sigma_i)=x^{\ast}\EE_0\]
is a fibre functor and the pair $(\St(T),\omega_x)$ is a neutral Tannakian category.
\end{prop}

Let us apply these notions to the given connected projective semi-stable curve $X$ with smooth generic fiber $X_K$, using the notation of the previous section.

\begin{defn} Let $\widebar K$ be a fixed algebraic closure of $K$,  $X_{\widebar K}=X_K\times_K \Sp(\widebar K)$ the base change and $x\in X_{\widebar K}(\widebar K)$ a closed point, we denote by $\pi^{\text{strat}}(X_{\widebar K},x)$ the affine group scheme associated with $(\St(X_{\widebar K}),\omega_x)$ via Tannakian duality.
\end{defn}

\begin{prop}[\cite{Kindtam}, Prop.2.15]\label{profinitecomplstrat} Let $\pi_1^{\text{\'et}}(X_{\widebar K},x)=\varprojlim_i \pi_i$ be the \'etale fundamental group of $X_{\widebar K}$, then there exists a morphism of $\widebar K$-group schemes  
\[\pi^{\textnormal{strat}}(X_{\widebar K},x)\to\varprojlim_i(\pi_i)_{\widebar K}=:\pi_1^{\text{\'et}}(X_{\widebar K},x)_{\widebar K}.\]
\end{prop}

We will now introduce the notion of stratified bundles with meromorphic descend data and generalize the results of the previous section to the category they form. 

\begin{defn} Given  $\mathcal{Y}$  a  geometric covering of $\widehat{X}$,  a coherent sheaf $\FF$ on $\mathcal{Y}$ is called \textit{meromorphic bundle} if there exists a locally free sheaf $\EE$ on $\mathcal{Y}$ such that  $\FF \otimes_A K\cong \EE\otimes_A K$.
\end{defn} 

\begin{rmk}\label{extfrob}
Note that, if $X$ is a projective semi-stable curve  over $S$  with geometrically connected smooth generic fibre and connected closed fibre, then so are its Frobenius twists $X^{(i)}$. Indeed, by \cite[Prop. 10.3.15.(a)]{Liu}, $X^{(i)}$ is a  projective semi-stable curve over $S^{(i)}$.
Moreover, the generic fibre of $X^{(i)}$ is   $(X^{(i)})_K\cong(X_K)^{(i)}$, which is clearly smooth and geometrically connected, and the closed fibre of $X^{(i)}$ is $(X^{(i)})_0\cong(X_0)^{(i)}$.
\end{rmk}

\begin{defn}  Given $(V,\rho)\in \Rp^{\textnormal{cts}}_K(\pi^{\textnormal{pro\'et}}_1(X_0,\xi))$, let $\mathcal{Y}_{\rho}^{(i)}$ and $\mathcal{Z}_{\rho}^{(i)}$ be the $i$-th Frobenius twists of $\mathcal{Y}_{\rho}$ and $\mathcal{Z}_{\rho}$ and $F_{Y/S}^i\colon \mathcal{Y}_{\rho}^{(i+1)}\to\mathcal{Y}_{\rho}^{(i)}$ the relative Frobenius over $S$.  A \textit{stratified bundle with meromorphic descent data on $\mathcal{Y}_{\rho}$} is given by the following data: 
\begin{itemize}
\item $\{\EE_i,h_w^i\}_{w\in \Aut(\mathcal{Y}_{\rho}|\widehat{X})}$,   meromorphic bundles on $\mathcal{Y}_{\rho}^{(i)}$ with meromorphic descent data 
\[h^i_w\colon\EE_i\otimes_A K\to w^{\star}\EE_i\otimes_A K,\] 
\item $\sigma_i$, isomorphisms of meromorphic descent data
\[\sigma_i \colon \{{F_{Y/S}^i}^{\star}\EE_{i+1},{F_{Y/S}^i}^{\star}h_w^{i+1}\}_{w\in \Aut(\mathcal{Y}_{\rho}|\widehat{X})}\to \{\EE_i,h_w^i\}_{w\in \Aut(\mathcal{Y}_{\rho}|\mathcal{Z}_{\rho})},\]
\end{itemize}
for each $i\geq 0$. 

In order to simplify the notation, we will often not specify the isomorphisms $\sigma_i$ and we will denote a stratified bundle with meromorphic descent data by $E=\{\EE_i,h_w^i\}$.
\end{defn}

\begin{defn}\label{mapstratmerom} A \textit{morphism of stratified bundles with meromorphic descent data} from $\{\EE_i,h_w^i,\sigma_i\}$ to $\{\GG_i,k_w^i,\tau_i\}$ is given by a sequence $\{\alpha_i\}$ of morphisms of sheaves with meromorphic descent data on $\mathcal{Y}^{(i)}_{\rho}$  such that 
the following diagram is commutative
\[
\begin{CD}
{F_{Y/S}^i}^{\star}\{\EE_{i+1},h_w^{i+1}\} @>{{F_{Y/S}^i}^{\star}\alpha_{i+1}}>>{F_{Y/S}^i}^{\star}\{\GG_{i+1},k_w^{i+1}\}\\
@V\sigma_iVV @VV\tau_iV \\
\{\EE_i,h_w^i,\sigma_i\} @>{\alpha_i}>> \{\GG_i,k_w^i,\tau_i\} \text{ .}
\end{CD}
\]

We denote by $\St^m(\mathcal{Y}_{\rho})$ the category of stratified bundle with meromorphic descent data on $\mathcal{Y}_{\rho}$.
\end{defn}

\begin{lem}\label{repstrat} A representation $(V,\rho)\in \Rp^{\textnormal{cts}}_K(\pi^{\text{pro\'et}}_1(X_0,\xi))$  induces a stratified bundle with meromorphic descent data  on ${\mathcal{Y}_{\rho}}$. 
\end{lem}
\proof Given a representation $(V,\rho)\in \Rp^{\textnormal{cts}}_K(\pi^{\text{pro\'et}}_1(X_0,\xi))$ with $V$ a $K$-vector space of rank $n$, we set
\[\gamma\colon\text{Aut}(\mathcal{Y}_{\rho}|\widehat{X})\to\mathbb{Z}^{\star r}\star G_1\star\cdots\star G_N\]
to be the composition of the isomorphism in \hyperref[autY]{Equation \ref*{autY}} and the inversion. We set then $\tilde{\rho}:=\rho\circ \gamma$ and we fix a base $V\simeq K^n$. We define the sheaf $\{\OO_{\mathcal{Y}_{\rho}^{(i)}}^n,h_w^{\rho,i}\}_{w\in\Aut(\mathcal{Y}_{\rho}|\widehat{X})}$  with meromorphic descent data, where $h_w^{\rho,i}$ are given by 
\[
\begin{CD}
\OO_{\mathcal{Y}_{\rho}^{(i)}}\otimes_A K^n\simeq\OO_{\mathcal{Y}_{\rho}^{(i)}}\otimes_A V@>{h^{\rho,i}_w}>>\OO_{\mathcal{Y}_{\rho}^{(i)}}\otimes_A V\simeq\OO_{\mathcal{Y}_{\rho}^{(i)}}\otimes_A K^n\\
 f\otimes v@>>>f\otimes\tilde{\rho}(w)(v)\text{ .}
\end{CD}
\]
By construction, it is clear that
\[{F_{Y/S}^i}^{\star}\{\OO_{\mathcal{Y}_{\rho}^{(i+1)}}^n,h_w^{\rho,i+1}\}=\{\OO_{\mathcal{Y}_{\rho}^{(i)}}^n,h_w^{\rho,i}\}.\] 
Hence, the sequence $\{\OO_{\mathcal{Y}_{\rho}^{(i)}}^n,h_w^{\rho,i}\}$ is a stratified bundle with meromorphic descent data  on ${\mathcal{Y}_{\rho}}$. 
\endproof 

\begin{defn}\label{meromstrat} A \textit{meromorphic stratified bundle} on $\widehat{X}$ is a sequence $\{{\GG}_i, \sigma_i\}$ of meromorphic bundles ${\GG}_i$ on $\widehat{X}^{(i)}$ and isomorphisms 
\[\sigma_i \colon {F_{\widehat{X}}^i}^{\star}{\GG}_{i+1}\otimes_A K\to {\GG}_i\otimes_A K.\]
A \textit{morphism of meromorphic stratified bundles} from $\{ {\GG}_i, \sigma_i\}$ to $\{ {\GG}'_i, \tau_i\}$ is given by a sequence $\{\varphi_i\}$ of morphisms
\[\varphi_i\colon {\GG}_i\otimes_A K\to {\GG}'_i\otimes_A K\]
that are compatible with $\sigma_i$ and $\tau_i$. We denote by $\St^K(\widehat{X})$ the category of meromorphic stratified bundle on $\widehat{X}$.
\end{defn}

\begin{prop}\label{descstrat}Given $(V,\rho)\in \Rp^{\textnormal{cts}}_K(\pi^{\text{pro\'et}}_1(X_0))$, let $q_X\colon \mathcal{Y}_{\rho}\to\widehat{X}$ be the constructed geometric covering, then the pullback $q_X^{\star}$ induces an equivalence of categories between the category $\text{Strat}^K(\widehat{X})$ of meromorphic stratified bundles and the category $\text{Strat}^m(\mathcal{Y}_{\rho})$ of stratified bundles on $\mathcal{Y}_{\rho}$ with meromorphic descent data.

In particular, the category $\text{Strat}^m(\mathcal{Y}_{\rho})$ is equivalent to the category of F-divided sheaves on $X_K$, which will be denoted by $\text{Fdiv}(X_K)$. 
\end{prop}
\proof 
We first prove that $q_X^{\star}$ is essentially surjective. 

Given $\{\EE_i,h_w^i, \sigma_i\}\in\St^m(\mathcal{Y}_{\rho})$, by \hyperref[findesc]{Theorem \ref*{findesc}}, for every $i$, the sheaf $\{\EE_i,h_w^i\}_{w\in \Aut(\mathcal{Y}_{\rho}|\widehat{X})}$ with meromorphic descent data descends to a coherent sheaf ${\GG}_i$ on $\widehat{X}$
By fpqc descent, ${\GG}_i$ are meromorphic bundles. 
Let $F_{\widehat{X}}$, $F_{\mathcal{Y}}$ be the relative Frobenii on $\widehat{X}$ and $\mathcal{Y}_{\rho}$ respectively and $q_X^{i}\colon\mathcal{Y}^{i}\to\widehat{X}^{(i)}$, then we have that
\[\textit{Hom}_{\OO_{\widehat{X}^{(i)}}}({F_{\widehat{X}}^i}^{\star}{\GG}_{i+1}\otimes_A K,{\GG}_i\otimes_A K)\simeq\textit{Hom}_{\OO_{\mathcal{Y}_{\rho}^{(i)}}}({q_X^{i}}^{\star}({F_{\widehat{X}}^i}^{\star}{\GG}_{i+1})\otimes_A K,{q_X^{i}}^{\star}{\GG}_i\otimes_A K).\]
Since   $F_{\widehat{X}}^i\circ q_X^{i}=q_X^{i+1}\circ F_{\mathcal{Y}}^i$,
and  ${q_X^{i}}^{\star}{\GG}_i\simeq\{ \EE_i,h_w^{i}\}$, we get
\[\textit{Hom}_{\OO_{\widehat{X}^{(i)}}}({F_{\widehat{X}}^i}^{\star}{\GG}_{i+1}\otimes_A K,{\GG}_i\otimes_A K)\simeq\textit{Hom}({F_{\mathcal{Y}}^i}^{\star}\{\EE_{i+1},h_w^{i+1}\},\{\EE_i,h_w^{i}\}).\]
Hence, $\sigma_i$ induces $\OO_{\widehat{X}}^{(i)}$-linear isomorphism 
$\varphi_i\colon F_{\widehat{X}}^{\star}\GG_{i+1}\otimes_A K\to\GG_i\otimes_A K$.
Moreover, by construction of $\varphi_i$, the isomorphism $q_X^{\star}{\GG}_i\otimes_A K\simeq\EE_i\otimes_A K$ makes the following diagram commute, for every $i$,
\[
\begin{CD}
{F_i}^{\star}q_X^{\star}{\GG}_{i+1} @>{}>>{F_i}^{\star}\EE_{i+1}\\
@V{q_X^{\star}\varphi_i}_iVV @VV\sigma_iV \\
q_X^{\star}{\GG}_i @>{}>> \EE_i\text{ .}\\
\end{CD}
\]
This implies that $\{q_X^{\star}\GG_i, q_X^{\star}\sigma_i\}$ and $\{\EE_i,h_w^i\}$ are isomorphic stratified bundles with meromorphic descent data.

Since $q_X$ is flat, clearly $q_X^{\star}$ is a faithful functor.
Let $\{\EE_i,\sigma_i\}$ and $\{\GG_i,\tau_i\}$ be two meromorphic stratified bundles  on $\widehat{X}$ and $\alpha_i\colon {q_X^{i}}^{\star}\EE_i\otimes_A K\to {q_X^{i}}^{\star}\GG_i\otimes_A K$ a morphism of stratified bundles with meromorphic descent data. Then, by \hyperref[findesc]{Theorem \ref*{findesc}}, there exists a corresponding morphism $\beta_i\colon\EE_i\otimes_A K\to\GG_i\otimes_A K$, for every $i$. In order to prove that $q_X^{\star}$ is full, it suffices to show that the following diagram commutes
\[
\begin{CD}
{F_{\widehat{X}}^i}^{\star}{\EE}_{i+1}\otimes_A K @>{F_{\widehat{X}}^i}^{\star}\beta_{i+1}>>{F_{\widehat{X}}^i}^{\star}\GG_{i+1}\otimes_A K\\
@V{\sigma_i}_iVV @VV\tau_iV \\
{\EE}_i\otimes_A K @>{\beta_i}>> \GG_i\otimes_A K\text{ .}\\
\end{CD}
\]
Since $F_{\widehat{X}}^i\circ q_X^{i}=q_X^{i+1}\circ F_{\mathcal{Y}}^i$, it is clear that ${q_X^{i}}^{\star}{F_{\widehat{X}}^i}^{\star}\beta_{i+1}$ corresponds to ${F_{\mathcal{Y}}^i}^{\star}\alpha_{i+1}$ via 
\[\textit{Hom}({q_X^{i}}^{\star}{F_{\widehat{X}}^i}^{\star}{\EE}_{i+1}\otimes_A K,{q_X^{i}}^{\star}\GG_i\otimes_A K)\simeq\textit{Hom}_{\OO_{\widehat{X}^{(i)}}}({F_{\widehat{X}}^i}^{\star}{\EE}_{i+1}\otimes_A K,{\GG}_i\otimes_A K).\] 
By hypothesis $\alpha_i$ is a a morphism of stratified bundles with meromorphic descent data and the following diagram commutes
\[
\begin{CD}
q_X^{\star}{F_{\widehat{X}}^i}^{\star}{\EE}_{i+1}\otimes_A K @>{F_{\mathcal{Y}}^i}^{\star}\alpha_{i+1}>>q_X^{\star}{F_{\widehat{X}}^i}^{\star}{\GG}_{i+1}\otimes_A K\\
@V{q_X^{\star}\sigma_i}_iVV @VVq_X^{\star}\tau_iV \\
q_X^{\star}{\EE}_i \otimes_A K@>{\alpha_i}>> q_X^{\star}\GG_i\otimes_A K\text{ .}\\
\end{CD}
\]
Thus, $\{\beta_i\}$ is a morphism of meromorphic stratified bundles.

Similarly, we can conclude, in analogy with \hyperref[equivalence]{Theorem \ref*{equivalence}}, that  the categories $\St^{K}(\widehat{X})$ and $\text{Fdiv}(X_K)$ are equivalent.
\endproof

\begin{prop} \label{functstr} The descent of stratified bundles with meromorphic descent data associated to continuous representations of $\pi_1^{\text{pro\'et}}(X_0,\xi)$ induces a tensor functor
\[\text{sp}_K\colon\Rp^{\textnormal{cts}}_K(\pi_1^{\text{pro\'et}}(X_0,\xi))\to\text{Fdiv}(X_K).\]
\end{prop}
\proof By \hyperref[descstrat]{Proposition \ref*{descstrat}}, given  $(V,\rho)\in \Rp^{\textnormal{cts}}_K(\pi^{\text{pro\'et}}_1(X_0))$, the stratified bundle with meromorphic descent data $\{\OO_{\mathcal{Y}_{\rho}^{(i)}}^n,h_w^{\rho,i}\}$ induced by $\rho$ on $\mathcal{Y}_{\rho}$ descends to a F-divided sheaf $\{\FF^i_{\rho}\}$ on $X_K$. Thus, we can define
\[\text{sp}_K(V,\rho):=\{\FF^i_{\rho}\}\in\text{Fdiv}(X_K).\]

Let  $\varphi\colon (V,\rho)\to(W,\tau)$ be a morphism of representations and assume that $\rho$ factors through the group $\mathbb{Z}^{\star r}\star G_1\star\cdots\star G_N$ and $\tau$ factors through the group $\mathbb{Z}^{\star r}\star H_1\star\cdots\star H_N$, then we denote by ${\mathcal{Y}_{\rho}}$ and  ${\mathcal{Y}_{\tau}}$ the   geometric coverings of $\widehat{X}$ associated with $\rho$ and $\tau$, as in \hyperref[yrho]{Definition \ref*{yrho}}. 
Moreover we set $G_i^{\rho,\tau}$ to be the image of the map $\pi_1^{\text{\'et}}(\widebar{C_i})\to G_i \times H_i$ and we associate with the $\pi_1^{\text{pro\'et}}(X_0,\xi))$-set $\mathbb{Z}^{\star r}\star G_1^{\rho,\tau}\star\cdots\star G_N^{\rho,\tau}$
 a geometric covering of $\widehat{X}$, which we call $\mathcal{Y}_{\rho,\tau}$.

We set 
\[\rho'\colon \mathbb{Z}^{\star r}\star G_1^{\rho,\tau}\star\cdots\star G_N^{\rho,\tau}\to\text{Aut}(V)\]
 to be the unique group morphism  such that $\rho'(w)=\rho(w)$ for every $w\in \mathbb{Z}^{\star r}$,
and  $\rho'(g_i,h_i)=\rho(g_i)$ for every $(g_i,h_i)\in G_i^{\rho,\tau}$ and every $i=1,\dots,N$. Similarly, we define $\tau'$. 
By construction, there exist maps 
\[p_{\rho}\colon \mathcal{Y}_{\rho,\tau}\to \mathcal{Y}_{\rho}\text{ and }p_{\tau}\colon \mathcal{Y}_{\rho,\tau}\to \mathcal{Y}_{\tau}.\]
and we have that
\[p_{\rho}^{\star}\{\OO_{{\mathcal{Y}}_{\rho}}^n,h_w^{\rho}\}=\{\OO^n_{{\mathcal{Y}}_{\rho,\tau}},h_w^{\rho'}\}\text{ and }p_{\tau}^{\star}\{\OO_{\mathcal{Y}_{\tau}}^m,h_w^{\tau}\}=\{\OO_{\mathcal{Y}_{\rho,\tau}}^m,h_w^{\tau'}\}.\] 

By \hyperref[findesc]{Theorem \ref*{findesc}}, for every $i$, we set $F_i(\varphi)$ to be the map corresponding to the morphism of meromorphic descent data
\[\alpha_{\varphi}\colon \{\OO_{\mathcal{Y}_{\rho,\tau}^{(i)}}^n,h_w^{\rho',i}\}_{w\in\text{Aut}(\mathcal{Y}_{\rho,\tau}|\widehat{X})}\to\{\OO_{\mathcal{Y}_{\rho,\tau}^{(i)}}^m,h^{\tau',i}_w\}_{w\in\text{Aut}(\mathcal{Y}_{\rho,\tau}|\widehat{X})}\] 
defined as follows:
\[
\begin{CD}
\OO_{\mathcal{Y}_{\rho}^{(i)}}\otimes_A V@>{\alpha_{\varphi}}>>\OO_{\mathcal{Y}_{\rho}^{(i)}}\otimes_A W\\
 f\otimes v@>>>f\otimes\varphi(v)\text{ .}
\end{CD}
\]
By construction, it is clear that  the collection $\{F_i(\varphi)\}$ induces a morphism of F-divided sheaves from $\{\FF^i_{\rho}\}$ to $\{\FF^i_{\tau}\}$. It remains to show that the functor we constructed is a tensor functor. 

Given $(V,\rho),(W,\tau)$ two continuous representations, let $\mathcal{Y}_{\rho,\tau}$ be the geometric covering defined   above. Then we define the representation
\[\rho'\otimes\tau'\colon\mathbb{Z}^{\star r}\star G_1^{\rho,\tau}\star\cdots\star G_N^{\rho,\tau}\to\text{Aut}(V\otimes W)\] 
and we  associate with it the stratified bundle on $\mathcal{Y}_{\rho,\tau}$ with meromorphic descent data $\{\OO_{\mathcal{Y}_{\rho,\tau}^{(i)}}^n,h_w^{\rho'\otimes\tau',i}\}$.
The tensor product of stratified bundles with meromorphic descent data is defined as follows
\[\{\OO_{\mathcal{Y}_{\rho,\tau}^{(i)}}^n,h_w^{\rho',i}\}\otimes\{\OO_{\mathcal{Y}_{\rho,\tau}^{(i)}}^m,h_w^{\tau',i}\}:=\{\OO_{\mathcal{Y}_{\rho,\tau}^{(i)}}^{nm},h_w^{\rho',i}\otimes h_w^{\tau',i}\},\]
hence it is clear that 
\[\{\OO_{\mathcal{Y}_{\rho,\tau}^{(i)}}^n,h_w^{\rho'}\}\otimes\{\OO_{\mathcal{Y}_{\rho,\tau}^{(i)}}^m,h_w^{\tau'}\}\simeq\{\OO_{\mathcal{Y}_{\rho\otimes\tau}^{(i)}}^{nm}, h_w^{\rho'\otimes\tau'}\}.\]

By construction, $\{\OO_{\mathcal{Y}_{\rho\otimes\tau}^{(i)}}^n,h_w^{\rho'}\}$ descends to $\text{sp}_K(\rho)$, $\{\OO_{\mathcal{Y}_{\rho,\tau}^{(i)}}^m,h_w^{\tau'}\}$ descends to $\text{sp}_K(\tau)$ and $\{\OO_{\mathcal{Y}_{\rho\otimes\tau}^{(i)}}^{nm}, h_w^{\rho'\otimes\tau'}\}$ descends to $\text{sp}_K(\rho\otimes\tau)$. Thus, it follows that 
\[\text{sp}_K(\rho)\otimes\text{sp}_K(\tau)\simeq\text{sp}_K(\rho\otimes\tau).\]
All the properties of tensor functor can be easily checked in a similar way.
\endproof

\begin{thm} \label{groupmor} For every finite extension $L$ of $K$, the functor $\text{sp}_K$ can be extended to a tensor functor
\[\text{sp}_L\colon\Rp^{\textnormal{cts}}_L(\pi_1^{\text{pro\'et}}(X_0,\xi))\to\text{Fdiv}(X_L),\]
 where $X_L=X_K\times_K \Sp(L)$.
Moreover, fixing $x\in X_{\widebar K}(\widebar K)$, it induces, up to canonical natural transformation, a morphism of group schemes
\[\text{sp}\colon\pi^{\text{strat}}(X_{\widebar K},x)\to (\pi_1^{\text{pro\'et}}(X_0,\xi))^{\textnormal{cts}}.\]
\end{thm}
\proof Let $(V,\rho)\in\Rp_{\widebar K}(\pi_1^{\text{pro\'et}}(X_0))$, then $\rho$ factors through the group $\mathbb{Z}^{\star r}\star G_1\star\cdots\star G_N$, which is finitely generated. Hence, there exists a finite field extension $K\subset L$ and $(V_L,\rho_L)\in\Rp_L(\pi_1^{\text{pro\'et}}(X_0))$  such that
\[(V_L,\rho_L)\otimes_L \widebar K=(V,\rho).\]
We set $A_L$ to be the integral closure of $A$ in $L$, $S_L=\Sp(A_L)$ and we set $X_{S_L}=X\times_S S_L$. By definition, $A_L$ is a complete discrete valuation ring, whose residue field is $k$ and whose fraction field is $L$. 
By base change, $X_{S_L}$ is a projective semi-stable curve with geometrically connected smooth generic fibre $X_L$ and connected closed fibre $X_0$.

We can apply \hyperref[functstr]{Proposition \ref*{functstr}} to $X_{S_L}$ and we can define a tensor functor $\text{sp}_L$ that associates to $(V_L,\rho_L)$ a F-divided sheaf on $X_L$.
Let $\text{bs}_L\colon X_{\widebar K}\to X_L$ be the base change, then by \hyperref[Katzthm]{Theorem \ref*{Katzthm}} we get a tensor functor $\text{sp}$ defined by
\[\text{sp}(V,\rho):=\text{bs}_L^{\ast}(\text{sp}_L(V_L,\rho_L))\in\St(X_{\widebar K}),\]
where $\St(X_{\widebar K})$ is the Tannakian category of stratified bundles over $X_{\widebar K}$.

We fix $x\in X_{\widebar K}(\widebar K)$ and we set $\omega_{x}$ to be the associated fibre functor of $\St(X_{\widebar K})$ and $\omega_{\pi}$ the fiber functor of $\Rp^{\textnormal{cts}}_{\widebar K}(\pi_1^{\textnormal{pro\'et}}(X_0,\xi))$ given by the forgetful functor. 
Since $X$ is proper and flat over $S$, there exists a specialization $\xi\in X_0$ of the $\widebar K$-point $x\in X$.
Given $\FF_{\rho}=\text{sp}(\widebar {K}^n,\rho)$. Hence, 
\[x^{\star}\FF_{\rho}=(\FF_{\rho,\xi}\otimes_{\OO_{X,\xi}} K \otimes_A \widebar K)_x.\]
Since the morphism $q_X\colon \mathcal{Y}\to\widehat{X}\to X$ is flat, there exists a local trivialization $\FF_{\rho,\xi}\otimes_A K\simeq \OO^n_{X,\xi} \otimes_A K$, which, by tensoring with $\OO^n_{\mathcal{Y},y}$ over $\OO_{X,\xi}$ for $y\in \mathcal{Y}$ such that $q_X(y)=\xi$, corresponds to the isomorphism $(q_X^{\star}\FF_{\rho}\otimes _A K)_y\simeq \OO_{\mathcal{Y},y}\otimes K^n$ defining $\FF_{\rho}$. Thus, we have
\[x^{\star}\FF_{\rho}\simeq  (\OO_{X,\xi}^n\otimes_A K\otimes_{\OO_{X,\xi}} \widebar K)_x \simeq {\widebar K}^n.\]
It remains to show that this isomorphism is functorial.

Given a morphism $\varphi\colon(V,\rho)\to(W,\tau)$  of $\widebar K$-linear representations, we set $\FF_{\rho}=\text{sp}(V,\rho)$ and $\GG_{\tau}=\text{sp}(W,\tau)$. To prove the functioriality of the above isomorphism it suffices to show that the following diagram commutes 
\[
\begin{CD}
\FF_{\rho,\xi} \otimes_A K@>>>\GG_{\tau,\xi} \otimes_A K\\
@VVV @VVV \\
\OO_{X,\xi}^n\otimes_A K@>>> \OO_{X,\xi}^m\otimes_A K{ .}\\
\end{CD}\] 
By descent, it suffices to show that the following diagram commutes
\[
\begin{CD}
q_X^{\star}\FF_{\rho}\otimes_A K@>>>q_X^{\star}\GG_{\tau}\otimes_A K\\
@VVV @VVV \\
\OO_{\mathcal{Y},y}^n\otimes_A K@>>> \OO_{\mathcal{Y},y}^m\otimes_A K\\
\end{CD}\] 
on a small neighborhood of $y$, which is true by construction of $\FF_{\rho}$ and $\GG_{\tau}$.

We conclude that there exists a natural isomorphism $\gamma$
\[\gamma\colon\omega_x \circ \text{sp}\simeq\omega_{\pi}.\]
Let $\omega_{\pi}'\colon=\gamma(\omega_{\pi})$, by \cite[Cor. 2.9]{Miln}, the functor $\text{sp}$ corresponds to  a morphism of group schemes 
\[\text{sp}\colon \pi^{\text{strat}}(X_{\widebar K})\to \pi(\Rp^{\textnormal{cts}}_K(\pi^{\text{pro\'et}}_1(X_0)),\omega_{\pi}').\]
Moreover, $\omega_{\pi}'$ and $\omega_{\pi}$ are naturally isomorphic, so we have that
\[\pi(\Rp^{\textnormal{cts}}_K(\pi^{\text{pro\'et}}_1(X_0)),\omega_{\pi}')\simeq (\pi_1^{\text{pro\'et}}(X_0))^{\text{cts}}\]
and, composing with this isomorphism, we get a morphism of group schemes
\[\text{sp}\colon\pi^{\text{strat}}(X_{\widebar K})\to (\pi_1^{\text{pro\'et}}(X_0))^{\text{cts}}.\] 
\endproof

\section{Compatibility with the \'etale specialization map}

Given $x\in X_{\widebar{K}}(\widebar K)$, we  denote by $\text{sp}_{SGA1}$ the specialization map constructed by Groethendieck in \cite{SGA1}
\[\text{sp}_{SGA1}\colon\pi_1^{\text{\'et}}(X_{\widebar{K}},x)\to\pi_1^{\text{\'et}}(X_0,\xi).\]
This  specialization morphism induces a functor
\[\text{sp}_{SGA1}\colon \Rp^{\text{cts}}_{\widebar{K}}(\pi_1^{\text{\'et}}(X_0,\xi))\to\Rp^{\text{cts}}_{\widebar{K}}(\pi_1^{\text{\'et}}(X_{\widebar{K}},x)).\]
Furthermore, by \hyperref[complet]{Proposition \ref*{complet}}, the pro-finite completion induces a fully faithful functor 
\[c\colon  \Rp^{\text{cts}}_{\widebar{K}}(\pi_1^{\text{\'et}}(X_0,\xi))\to  \Rp^{\text{cts}}_{\widebar{K}}(\pi_1^{\text{pro\'et}}(X_0,\xi)).\]

Let $(\widebar{K}^n,\rho)\in\Rp^{\text{cts}}_{\widebar{K}}(\pi_1^{\text{\'et}}(X_{\widebar{K}},x))$, by continuity  and \hyperref[repdiscrete]{Lemma \ref*{repdiscrete}}, $\rho$ factors through a finite quotient $\pi_{\rho}$ of $\pi_1^{\text{\'et}}(X_{\widebar{K}},x)$. In particular, there exists a finite Galois cover $W_{\widebar{K}}$ of $X_{\widebar{K}}$ such that 
\[\text{Aut}(W_{\widebar{K}}|X_{\widebar{K}})=\pi_{\rho}^{\text{op}}.\]
We can define  descend data $\{h^{\rho}_g\}_{g\in \pi_{\rho}}$ for the sheaf $\OO_{W_{\widebar{K}}}^n$ on $W_{\widebar{K}}$ as follows
\[
\begin{CD}
\OO_{W_{\widebar{K}}}^n @>{h^{\rho}_g}>>\OO_{W_{\widebar{K}}} ^n  \\
 (f_i)@>>>\rho(g)(f_i)\text{ .}
\end{CD}
\] 
Since ${W_{\widebar{K}}}\to X_{\widebar{K}}$ is a morphism of effective descent for coherent sheaves,  $\{\OO_{W_{\widebar{K}}}^n,h^{\rho}_g\}$ descends to a coherent sheaf $\EE$ on $X_{\widebar{K}}$ that, by construction, is locally free.
As in  the proof of \hyperref[functstr]{Proposition \ref*{functstr}}, if we repeat the argument for the Frobenius twists of $X_{\widebar{K}}$ we can define a functor
\[F\colon \Rp^{\text{cts}}_{\widebar{K}}(\pi_1^{\text{\'et}}(X_{\widebar{K}},x))\to\St(X_{\widebar{K}}).\]

Combining these functors with the specialization functor, we get the following diagram:
\begin{equation}\label{commdiag}
\begin{CD}
\Rp^{\textnormal{cts}}_{\widebar{K}}(\pi_1^{\text{\'et}}(X_0,\xi)) @>\text{sp}_{SGA1}>>\Rp^{\text{cts}}_{\widebar{K}}(\pi_1^{\text{\'et}}(X_{\widebar{K}},\epsilon))\\
@VV{c}V @VVFV \\
\Rp^{\text{cts}}_{\widebar{K}}(\pi_1^{\text{pro\'et}}(X_0,\xi)) @>\text{sp}>>\St(X_{\widebar{K}})\text{ .}\\
\end{CD}
\end{equation}

\begin{lem}\label{pcommdiag} The diagram \eqref{commdiag} is commutative up to a natural transformation.
\end{lem}
\proof Let  $(V,\rho)\in\Rp^{\text{cts}}_{\widebar{K}}(\pi_1^{\text{\'et}}(X_0,\xi))$, then by continuity and \hyperref[repdiscrete]{Lemma \ref*{repdiscrete}}, the morphism $\rho$ factors through a finite quotient $G_{\rho}$ of $\pi_1^{\text{\'et}}(X_{0},\xi)$. Moreover, there exists a finite field extension $L$ of $K$ and $(V_L,\rho_L)\in\Rp^{\text{cts}}_{L}(\pi_1^{\text{\'et}}(X_0,\xi))$ such that
\[(V,\rho)=(V_L,\rho_L)\otimes_K \widebar{K}.\]
For simplicity we  call $\rho$ also the representation with coefficients in $L$.

Hence, we have the following commutative diagram
\[\begin{tikzpicture}[node distance=2.5cm, auto]
\node (P) {$\pi_1^{\text{pro\'et}}(X_0,\xi)$};
\node(Q)[right of=P] {$\pi_1^{\text{\'et}}(X_0,\xi)$};
\node (B) [below of=P] {$G_{\rho}$};
\node (C) [right of=B] {$\text{Aut}(V_L)\text{ .}$};
\draw[->](P) to node {$c$}(Q); 
\draw[->](Q) to node {$$}(B);
\draw[->](Q) to node {$\rho$}(C);
\draw[->] (P) to node {$p$} (B);
\draw[->] (B) to node {$\bar{\rho}$} (C);
\end{tikzpicture}
\]

Moreover, the morphism $p$ factors through the quotient $\mathbb{Z}^{\star r}\star H_1\star\cdots\star H_N$,
where, if $j$ is the natural morphism $j\colon\pi_1^{\text{\'et}}(\widebar{C_j})\to\pi_1^{\text{pro\'et}}(X_0,\xi)$ and $p_j=p\circ j$, 
\[H_j=\pi_1^{\text{\'et}}(\widebar{C_j})/p_j^{-1}(\text{Id}).\]

We recall that to define  $\text{sp}(\widebar K, \rho\circ c)$,  we have set
\[G_j=\pi_1^{\text{\'et}}(\widebar{C_j})/(\rho\circ c \circ j)^{-1}(\text{Id}).\]
Since $\rho\circ c \circ j=\bar{\rho}\circ p \circ j$ and  $\bar{\rho}$ is injective, we have that $H_j=G_j$. Thus, there exists a $\pi_1^{\text{pro\'et}}(X_0,\xi)$-equivariant morphism 
\[q\colon \mathbb{Z}^{\star r}\star H_1\star\cdots\star H_N\to G_{\rho}\] that, together with the quotient map $\pi_1^{\text{pro\'et}}(X_0,\xi)\to  \mathbb{Z}^{\star r}\star H_1\star\cdots\star H_N$ completes the above commutative diagram. 

Let $X_{S_L}$ be defined as \hyperref[groupmor]{Theorem \ref*{groupmor}} and let $\widehat{X}_{S_L}$ be the completion of $X_{S_L}$ along $X_0$, then the  set $\mathbb{Z}^{\star r}\star G_1\star\cdots\star G_N$ corresponds to $\mathcal{Y}_{\rho}$, while $G_{\rho}$ correspond to a geometric coverings of $\widehat{X}_{S_L}$ that we call $\mathcal{W}$. Moreover, $q$ corresponds to a  $\widehat{X}_{S_L}$-morphism 
\[
\begin{tikzcd}
 & \mathcal{W} \arrow{dr}{p_{\mathcal{W}}} \\
\mathcal{Y}_{\rho}\arrow{ur}{q} \arrow{rr}{p_{\mathcal{Y}}} && \widehat{X}_{S_L} .
\end{tikzcd}
\] 

By construction,  $\text{sp}_L(V_L,\rho\circ c)$ corresponds to a sequence $\{\FF_{\rho}^i\}$ of meromorphic bundles on $\widehat{X}_{S_L}$ such that
\[p_{\mathcal{Y}_{\rho}}^{\star}\{\FF_{\rho}^i\}\simeq\{\OO^n_{\mathcal{Y}_{\rho}^{(i)}},h_w^{\bar{\rho}\circ q,i}\}_{w\in(\mathbb{Z}^{\star r}\star G_1\star\cdots\star G_N)^{\text{op}}}.\]

If ${W}$ is the finite \'etale covering of $X_{S_L}$ corresponding to $\mathcal{W}$ and $W_{\widebar K}$ is its geometric generic fibre, since ${W}$ is normal, by \cite[Lemma 4.11]{Lorenz}, we deduce that 
\[\text{Aut}(W_{\widebar K}|X_{\widebar K})\simeq G_{\rho}^{\text{op}}.\]
By definition of the functor $F$ and \hyperref[equivalence]{Theorem \ref*{equivalence}}, $F(\text{sp}_{SGA1}(V_L,\rho))$ corresponds to a sequence of meromorphic bundles $\{\GG^{\rho}_i\}$ on $\widehat{X}_{S_L}$ such that
\[p_{\mathcal{W}}^{\star}\GG_i^{\rho}\simeq\{\OO^n_{\mathcal{W}^{(i)}},h_g^{\bar{\rho},i}\}_{g\in G_{\rho}^{\text{op}}}.\]

It is  easy to see that $\{\OO^n_{\mathcal{Y}_{\rho}^{(i)}},h_w^{\bar{\rho}\circ q,i}\}$ descends to $\{\OO^n_{\mathcal{W}^{(i)}},h_g^{\bar{\rho},i}\}$ on $\mathcal{W}$. Indeed, we have $\Aut(\mathcal{Y}_{\rho}|\mathcal{W})=\text{ker}(q)^{\text{op}}$ and 
\[\{\OO^n_{\mathcal{Y}_{\rho}^{(i)}},h_w^{\bar{\rho}\circ q,i}\}_{w\in\text{ker}(q)^{\text{op}}}=\{\OO^n_{\mathcal{Y}_{\rho}^{(i)}},\text{Id}\}_{w\in\text{ker}(q)^{\text{op}}},\]
which implies that $\{\OO^n_{\mathcal{Y}_{\rho}^{(i)}},h_w^{\rho,i}\}_{w\in\text{ker}(q)^{\text{op}}}$ descends to the trivial stratified bundle $\{\OO_{\mathcal{W}^{(i)}}^n\}$.
Moreover, for every $w\in(\mathbb{Z}^{\star r}\star G_1\star\cdots\star G_N)^{\text{op}}$ such that $q(w)=g$, $h_w^{q\circ\bar{\rho},i}$ corresponds to $h_g^{\bar{\rho},i}$ via the identification
\[\textit{Hom}(\{\OO^n_{\mathcal{Y}_{\rho}^{(i)}},\text{Id}\},\{\OO^n_{\mathcal{Y}_{\rho}^{(i)}},\text{Id}\})\simeq\textit{Hom}(\{\OO^n_{\mathcal{W}^{(i)}}\},\{\OO^n_{\mathcal{W}^{(i)}}\}).\]

Hence, by construction of the functors $F$ and \text{sp}, we find that
\[ \text{sp}(c(V_L,\rho))= F(\text{sp}_{SGA1}(V_L,\rho)).\]
\endproof

\begin{prop} The diagram \eqref{commdiag} induces, up to conjugation by a rational point, the following commutative diagram of group schemes 
\begin{equation}
\begin{CD}
\pi^{\text{strat}}(X_{\widebar{K}},x) @>F>>\pi_1^{\text{\'et}}(X_{\widebar{K}},\epsilon)_{\widebar{K}}\\
@VV{\text{sp}}V @VV\text{sp}_{SGA1}V \\
\pi_1^{\text{pro\'et}}(X_0,\xi)^{\textnormal{cts}} @>c>>\pi_1^{\text{\'et}}(X_0,\xi)_{\widebar{K}},\\
\end{CD}
\end{equation}
where $\pi_1^{\text{\'et}}(X_{\widebar{K}},\epsilon)_{\widebar{K}}$ and $\pi_1^{\text{\'et}}(X_0,\xi)_{\widebar{K}}$ are defined as in \hyperref[commlim]{Lemma \ref*{commlim}}.
\end{prop}
\proof Since the construction of the functor  $F$  is analogous to the construction of the specialization functor, following the same reasoning of \hyperref[groupmor]{Theorem \ref*{groupmor}}, we can conclude that  $F$ induces a morphism between the corresponding group schemes.
Then the commutativity of the above diagram of group schemes  follows from the previous proposition. 
\endproof

\section*{Bibliography}
\bibliographystyle{elsarticle-num}
\bibliography{biblio}

\end{document}